\documentclass[10pt,twoside]{amsart}
\usepackage{amsmath}
\usepackage{amssymb}
\usepackage{mathrsfs}		
\usepackage{graphicx}
\usepackage[curve]{xypic}
\usepackage{leftidx}		
\makeatletter
\def\@tocline#1#2#3#4#5#6#7{\relax
\ifnum #1>\c@tocdepth 
  \else 
    \par \addpenalty\@secpenalty\addvspace{#2}%
\begingroup \hyphenpenalty\@M
    \@ifempty{#4}{%
      \@tempdima\csname r@tocindent\number#1\endcsname\relax
 }{%
   \@tempdima#4\relax
 }%
 \parindent\z@ \leftskip#3\relax \advance\leftskip\@tempdima\relax
 \rightskip\@pnumwidth plus4em \parfillskip-\@pnumwidth
 #5\leavevmode\hskip-\@tempdima #6\nobreak\relax
 \ifnum#1<0\hfill\else\dotfill\fi\hbox to\@pnumwidth{\@tocpagenum{#7}}\par
 \nobreak
 \endgroup
  \fi}
\makeatother

\usepackage[colorlinks=true, urlcolor=blue,bookmarks=true,bookmarksopen=true, citecolor=blue]{hyperref}

\numberwithin{equation}{section}

\newtheorem{theorem}{Theorem}[section]
\newtheorem{defn}[theorem]{Definition}

\newtheorem{corollary}[theorem]{Corollary}

\newtheorem{example}[theorem]{Example}

\newtheorem{lemma}[theorem]{Lemma}
\newtheorem{prop}[theorem]{Proposition}
\newtheorem{remark}[theorem]{Remark}

\def \begineq{\begin{equation}}
\def \endeq{\end{equation}}

\def \bb{\mathbb}

\def \mc{\mathcal}
\def \mf{\mathfrak}

\newcommand{\Cnabla}[1]{{\leftidx{^C}{\nabla}{^{#1}}}}

\newcommand{\bmf}[1]{{\bar{\mf #1}}}

\def \CC{{\bb{C}}}

\def \GG{{\bb G}}
\def \HH{{\bb H}}
\def \JJ{{\bb{J}}}

\def \RR{{\bb{R}}}

\def \TT{{\bb{T}}}

\def \ZZ{{\bb{Z}}}

\def \LLC{{\mc L}}

\def \({\left(}
\def \){\right)}
\def \<{\left\langle}
\def \>{\right\rangle}

\def \bar{\overline}
\def \bdsum{\mathop{\boxplus}}

\def \deg{\mathrm{deg}}
\def \dsum{\mathop{\oplus}}
\def \Dsum{\mathop{\oplus}}

\def \inter{\cap}
\def \into{\hookrightarrow}

\def \tensor{\otimes}
\def \union{\cup}

\def \varleq{\leqslant}

\def \xto{\xrightarrow}

\def \Ad{{\rm Ad}}
\def \a{{a}}
\def \ad{{\rm ad}}

\def \Aut{{\rm Aut}}
\def \Cl{{\rm Cl}}

\def \End{{\rm End}}

\def \id{{id}}

\def \Pic{{\rm Pic}}
\def \PPic{{\mathbb P{\rm ic}}}

\def \Span{{\rm Span}}
\def \Spin{{\rm Spin}}

\def \tr{{\rm tr}}
\def \vol{{\rm vol}}

\renewcommand{\1}{1\!\!1}

\def \qed{\hfill $\square$ \vspace{0.02in}}

\begin{document}
\title{On generalized K\"ahler geometry on compact Lie groups}

\author{Shengda Hu}
\address{Department of Mathematics, Wilfrid Laurier University, 75 University Ave. West, Waterloo, Canada}
\email{shu@wlu.ca}

\abstract
 We present some fundamental facts about a class of generalized K\"ahler structures defined by invariant complex structures on compact Lie groups. The main computational tool is the BH-to-GK spectral sequences that relate the bi-Hermitian data to generalized geometry data. The relationship between generalized Hodge decomposition and generalized canonical bundles for generalized K\"ahler manifolds is also clarified.
\endabstract

\maketitle

\tableofcontents

\section{Introduction}\label{sect:intro}

Generalized K\"ahler structures were first introduced by Gualtieri in his thesis \cite{Gualtieri04}, as the analogue to classical K\"ahler structures in the framework of generalized geometry \`a la Hitchin \cite{Hitchin02}. Recall that a \emph{generalized complex structure} on a manifold $M$ is an integrable almost complex structure on $\TT M : = TM \dsum T^*M$, on which a structure of Courant algebroid can be defined by a closed three form $\gamma \in \Omega^3(M)$. A \emph{generalized K\"ahler structure} is a pair of generalized complex structures $(\JJ_+, \JJ_-)$, satisfying a compatibility condition. In \cite{Gualtieri04} it's shown that generalized K\"ahler structures are equivalent to  \emph{bi-Hermitian structures with torsion} as given by Gates-Hull-Ro\v cek \cite{GatesHullRocek}, where it was argued that they are the most general backgrounds for $N = (2, 2)$ supersymmetry.

Let $K$ be a compact Lie group of dimension $2n$, with Lie algebra $\mf k$. 
It is well-known \cite{Gualtieri04} that $K$ admits natural generalized K\"ahler structures. In terms of the corresponding bi-Hermitian structures, they are defined by left and right invariant complex structures on $K$, which exist by Samelson \cite{Samelson53} and Wang \cite{Wang54}. More precisely, let $I_+$ be a right invariant complex structure and $I_-$ a left invariant one, we have
$$\JJ_\pm := \frac{1}{2} \begin{pmatrix}I_+ \pm I_- & -(\omega_+^{-1} \mp \omega_-^{-1} \\ \omega_+ \mp \omega_- & -(I_+^* \pm I_-^*)\end{pmatrix} : \TT M \to \TT M$$
where $\omega_\pm$ are the K\"ahler forms of $I_\pm$ given by a bi-invariant metric $\sigma$ on $K$. These generalized K\"ahler structures are in general non-K\"ahler, while the metric $\sigma$ is always Gauduchon with respect to both complex structures $I_\pm$ (cf. Lemma \ref{lemma:alwaysGauduchon}). In particular, degrees are well defined for $I_\pm$-holomorphic vector bundles, as well as generalized holomorphic vector bundles on $K$ (Hu-Moraru-Seyyedali \cite{HuMoraruSeyyedali}).

In this note, we present some fundamental facts on a class of generalized K\"ahler structures on a compact Lie group $K$.
Let $\gamma \in\Omega^3(K)$ be the Cartan $3$-form. The corresponding Courant algebroid $\TT K$ can be trivialized by the natural extended actions of $K \times K$ (Alekseev-Bursztyn-Meinrenken \cite{AlekseevBursztynMeinrenken}, cf. Meinrenken \cite{Meinrenken13}). Corresponding to the invariant K\"ahler structures on an even dimensional torus $T$, we have the notion of \emph{Lie algebraic generalized K\"ahler structures} on $K$.
\begin{defn}\label{defn:recallLiealgebraicgenKahler}
 A generalized complex structure on $K$ defined via the Courant trivialization by a complex Lagrangian subalgebra $\mf L \subset \mf d$ satisfying $\mf L \inter (\mf k\dsum \mf k) = \{0\}$ is called a \emph{Lie algebraic generalized complex structure}. A generalized K\"ahler structure $(K, \gamma; \JJ_\pm, \sigma)$ is a \emph{Lie algebraic generalized K\"ahler structure} if both $\JJ_\pm$ are Lie algebraic generalized complex structures.
\end{defn}
\noindent
It turns out that such generalized K\"ahler structures are precisely those defined by the invariant complex structures. We comment that these structures are different from the \emph{invariant structures} considered in the literature, e.g. Alekseevsky-David \cite{AlekseevskyDavid10} and the references therein. In general, Lie algebraic generalized K\"ahler structures are not invariant under the left or right actions of the group; and a Lie algebraic generalized complex structure may not be part of a Lie algebraic generalized K\"ahler structure (cf. Lemma \ref{lemma:LiealgebraicgenKahler}).

Many classical notions have their generalized analogues, such as generalized Calabi-Yau manifold \cite{Hitchin02}, generalized holomorphic vector bundles, generalized canonical bundles \cite{Gualtieri04} and generalized Hodge decompositions (Gualtieri \cite{Gualtieri0409}, Cavalcanti \cite{Cavalcanti11}), to name a few. We work out some of these analogues for Lie algebraic generalized K\"ahler structures. For example, when $\pi_1(K)$ is torsionless, we show (Theorem \ref{thm:canonicalbundledegrees}) that the degrees of the canonical line bundles are invariants of the group $K$, (upto the choice of the bi-invariant metric $\sigma$).
\begin{theorem}\label{thm:degreecanonicalbundle}
 Let $\mf g = \mf k \tensor \CC$ be the complexified Lie algebra. With a choice of Cartan subalgebra $\mf h \subset \mf g$ and the corresponding set of positive roots $R_+$, let $\displaystyle{\rho = \frac{1}{2} \sum_{\alpha \in R_+} h_\alpha \in \mf h}$. The degree of the canonical line bundle for $\JJ_+$ is given by $-2|\rho|^2$, where the norm is induced by the bi-invariant metric on $K$.
\end{theorem}
\noindent
As a corollary (Corollary \ref{coro:canonicalgenKahlerisCalabiYau}), we recover the result in Cavalcanti \cite{Cavalcanti12}, that $\JJ_+$ is not generalized Calabi-Yau unless $K$ is a torus, in which case, it is part of a generalized Calabi-Yau metric structure (Corollary \ref{coro:torusgenCalabiYau}).

Let $L_+$ be the $i$-eigenbundle of the generalized complex structure $\JJ_+$, which decomposes further into $\pm i$-eigensubbundles of $\JJ_-$:
$$L_+ = \ell_+ \dsum \ell_-$$
then $\ell_\pm$ form a \emph{matched pair} of Lie algebroids (cf. Laurent-Gengoux-Sti\'enon-Xu \cite{LaurentGengouxStienonXu08}, Mackenzie \cite{Mackenzie07}, Mokri \cite{Mokri97} and Lu \cite{Lu97}). The double complex $(\Omega^{p,q}(\bar L_+), \bar \delta_+, \bar\delta_-)$ associated to the matched pair $(\bar \ell_+, \bar\ell_-)$ appeared first in Gualtieri \cite{Gualtieri10} in describing the deformation theory for $\JJ_\pm$. The associated spectral sequences are called the \emph{BH-to-GK spectral sequences} (\S \ref{sect:liealgebroidspectralsequence}). In Lemma \ref{lemma:explicitreducedholomorphicstructure}, we show that the differentials $\bar\delta_\pm$ can be explicitly identified using components of the corresponding Bismut connections.

The \emph{BH-to-GK} spectral sequence for $\JJ_+$-holomorphic vector bundles provides the relation between the generalized Hodge decomposition and cohomology of the generalized canonical bundles.
Recall  that for a generalized complex structure $\JJ$ on $M$ with $i$-eigenbundle $L$, the \emph{generalized canonical line bundle} $U \subset \wedge^* T_\CC^*M$ is generated locally by the pure spinors $\chi \in \Omega^*(M)$ defining $\bar L$. It is naturally a $\JJ$-holomorphic line bundle \cite{Gualtieri04}. Furthermore, suppose that $\JJ = \JJ_+$ is part of a generalized K\"ahler structure, then $U_{-n, 0} := C^\infty(U) \subset \Omega^*(M)$ generates, via the spinor actions of $\ell_\pm$, the \emph{generalized Hodge decomposition} \cite{Gualtieri0710} of the twisted de Rham complex $(\Omega^*(M), d_\gamma := d + \gamma \wedge)$. In particular, $d_\gamma$ decomposes (cf. \eqref{eq:twisteddiffdecompositions})
$$d_\gamma = \delta_+^\gamma + \delta_-^\gamma + \bar\delta_+^\gamma + \bar\delta_-^\gamma$$
For a generalized K\"ahler manifold, we have (cf. Theorem \ref{thm:spectralsequenceaspartofHodge})
\begin{theorem}\label{thm:introHodgedecomposition}
 The double complex inducing the BH-to-GK spectral sequence is naturally isomorphic to the \emph{half} of the generalized Hodge decomposition given by $(\bar\delta_+^\gamma, \bar\delta_-^\gamma)$. More precisely
 $$(\Omega^{p,q}(U; \bar L), D_+, D_-) \cong (U_{r,s}, \bar\delta_+^\gamma, \bar\delta_-^\gamma)$$
 where $r = p+q - n$ and $s = p-q$.
\end{theorem}
\noindent
When $K$ is semi-simple, it is well-known that the \emph{twisted de Rham cohomology} $H_\gamma^*(K)$ vanishes (e.g. Ferreira \cite{Ferreira13}). We then obtain the following vanishing result (cf. Proposition \ref{prop:canonicalcohotrivial}).
\begin{theorem}\label{thm:vanishingcanonicalcohomology}
 For a semi-simple Lie group $K$ with torsionless $\pi_1$, let $\JJ$ be part of a Lie algebraic generalized K\"ahler structure and $U$ the generalized canonical line bundle, then $H^*(U; \JJ) = 0$.
\end{theorem}
\noindent
For general compact Lie groups, the Hodge decomposition can be explicitly described at the level of Lie algebras, via the Courant trivialization (\S\ref{subsect:CliffordHodgedecomposition}).

In \cite{Gualtieri10}, it is shown that holomorphic reduction induces on $\bar \ell_\pm$ the structures of $I_\mp$-holomorphic Lie algebroids, which are denoted $\mc A_\mp$. Using Morita equivalence, \cite{Gualtieri10} further showed that the category of $\JJ_+$-holomorphic vector bundles is equivalent to the categories of (locally free) holomorphic $\mc A_\pm$-modules (cf. Corollary $3.15$, \emph{loc. cit.}). By the explicit identification of the differentials $\bar\delta_\pm$, we obtain the following refinement (cf. Proposition \ref{prop:isomorphiccategories}).
\begin{theorem}\label{thm:isomorphiccategories}
 On a generalized K\"ahler manifold, the categories of $\JJ_+$-holomorphic vector bundles and (locally free) $\mc A_\pm$-modules are isomorphic to each other.
\end{theorem}
\noindent
We note that the above theorem can be obtained as special case from results in \cite{LaurentGengouxStienonXu08}, while we provide details for completeness.

A key simplification for compact Lie groups comes from the observation that $\mc A_\mp$ are trivial as holomorphic vector bundles (cf. Proposition \ref{prop:Liealgebroidastrivialholobundle}). As further applications of the BH-to-GK spectral sequence, we compute the Lie algebroid cohomology $H^*(\bar L_+)$ and the $\JJ_+$-Picard group $\PPic^+_0(K)$ of $\JJ_+$-holomorphic line bundles (cf. Corollary \ref{coro:cohomologyofLonK} and Proposition \ref{prop:JPicardonK}).
\begin{theorem}\label{thm:cohoandpicard}
 Suppose that $\JJ = \JJ_+$ is part of a generalized K\"ahler structure on $K$, and let $L := L_+$ denote the $i$-eigenbundle of $\JJ$. Then 
 $$H^*(\bar L_+) \cong \wedge^* \CC^{2r} \text{ and } \PPic^+_0(K) \cong \Pic_0^+(K) \times \CC^r$$ 
 where $2r$ is the rank of $K$ and $\Pic^+_0(K)$ denotes the identity component of the Picard group of $I_+$-holomorphic line bundles on $K$.
\end{theorem}

The structure of the paper is as follows. We briefly review the general facts about generalized K\"ahler structures in \S \ref{sect:generalgenKahler} and recall the Courant trivialization as well as the differential in the Clifford algebra in \S \ref{sect:Couranttrivialization}. In \S \ref{sect:liealgebroidspectralsequence}, we describe the BH-to-GK spectral sequences and the relation to the generalized Hodge decomposition. The more detailed computations for compact Lie groups are contained in the last section \S \ref{sect:canonicalgenkalher}, where we also identify Lie algebraic generalized K\"ahler structures.

{\bf Acknowledgement:} I'd like to thank Kaiming Zhao for helpful discussions. I'd like to thank Marco Gualtieri, Ruxandra Moraru and Reza Seyyedali for collaborations that led to considerations in this article, as well as helpful discussions.

\section{Generalized K\"ahler geometry} \label{sect:generalgenKahler}
Let $M$ be a smooth manifold of dimension $2n$, $\gamma \in \Omega^3(M)$ a closed $3$-form, and $g$ a Riemannian metric on $M$. We give a self-contained description of basic facts on generalized K\"ahler geometry and generalized holomorphic vector bundles. More details can be found in several of Gualtieri's papers \cite{Gualtieri04}, \cite{Gualtieri0703}, \cite{Gualtieri0710} and \cite{Gualtieri10}.

\subsection{Equivalent descriptions} \label{subsect:genKahlergeomdescription}
We recall here two equivalent descriptions of generalized K\"ahler geometry: first as a \emph{bi-Hermitian structure} \cite{GatesHullRocek}; then as the analogue to the K\"ahler structures in the context of generalized geometry \cite{Gualtieri04}.

Let $I_\pm$ be integrable almost complex structures on $M$ such that $g$ is Hermitian with respect to both of them. Let $\omega_\pm = g\circ I_\pm$ be the corresponding K\"ahler form. Then the tuple $(M, \gamma; g, I_\pm)$ defines a \emph{bi-Hermitian structure} if
$$\pm d^c_\pm \omega_\pm = \gamma, \text{ where } d^c = i(\bar\partial - \partial)$$

For each Hermitian manifold, there is a corresponding \emph{Bismut connection}, for which both the metric and the integrable complex structure are parallel. Here, the corresponding Bismut connections are denoted $\nabla^\pm$, then
$$\nabla^\pm g = 0 \text{ and } \nabla^\pm I_\pm = 0$$
The relation to the Levi-Civita connection is given by $\gamma$, as will be shown below
$$\nabla^\pm_X Y = \nabla_X Y \pm \frac{1}{2}g^{-1}\iota_X\iota_Y \gamma$$

Consider the generalized tangent bundle $\TT M = TM \dsum T^*M$, which admits a Dorfman bracket defined by $\gamma$:
$$(X+\xi)*(Y+\eta) = [X,Y] + \LLC_X\eta - \iota_Y d\xi + \iota_X\iota_Y \gamma$$
The natural projection $\TT M \to TM$ is denoted $a$.
A \emph{generalized almost complex structure} $\JJ : \TT M \to \TT M$ satisfies $\JJ^2 = -\1$ and is orthogonal with respect to the natural pairing 
$$\<X+\xi, Y+\eta\> = \frac{1}{2}(\iota_X\eta + \iota_Y \xi)$$
Let $L$ denote the $i$-eigensubbundle of $\JJ$ in $\TT_\CC M = \TT M\tensor_\RR \CC$, then it is maximally isotropic with respect to the pairing $\<,\>$. The structure $\JJ$ is \emph{integrable} and is called a \emph{generalized complex structure} when $L$ is \emph{involutive} with respect to the Dorfman bracket. In other words, $L$ is a \emph{Dirac structure}.

A pair of generalized complex structures $(\JJ_+, \JJ_-)$ defines a \emph{generalized K\"ahler structure} (with corresponding metric $g$) if they commute and such that
$$\GG = - \JJ_+\JJ_- = -\JJ_-\JJ_+ = 
\begin{pmatrix}
 0 & g^{-1} \\ g & 0
\end{pmatrix}
$$
A generalized K\"ahler structure induces the following decomposition:
$$\TT_\CC M = C_+^\CC \dsum C_-^\CC = L_+ \dsum \bar L_+ = L_- \dsum \bar L_- = \ell_+ \dsum \ell_- \dsum \bar \ell_+ \dsum \bar \ell_-$$
where $C_\pm$ are the graphs of $\pm g$ and $C_\pm^\CC$ their complexification, $L_+$ and $L_-$ are respectively the $i$-eigensubbundles of $\JJ_+$ and $\JJ_-$, $\ell_+ = L_+ \inter L_-$ and $\ell_- = L_+ \inter \bar L_-$. We have $C_\pm^\CC = \ell_\pm \dsum \bar\ell_\pm$, which defines complex structures $I_\pm$ on $TM$ by the restriction of $\JJ_+$, i.e. for $X \in TM$
$$I_\pm(X) = a \left(\JJ_+ (X \pm g(X))\right)$$
Then $a$ gives isomorphisms of $\ell_\pm$ to the holomorphic tangent bundles $T^\pm_{1,0} M$
$$\ell_\pm = \{X \pm g(X) = X \mp i\iota_X\omega_\pm : X \in T^\pm_{1,0} M\}$$
The condition that $\ell_\pm$ are involutive is equivalent to
$$(X \mp i\iota_X \omega_\pm) * (Y \mp i\iota_Y \omega_\pm) = [X, Y] \mp i(\LLC_{X}\iota_Y\omega_\pm - \iota_{Y}d\iota_X \omega_\pm) + \iota_X\iota_Y \gamma = [X, Y] \mp i\iota_{[X,Y]} \omega_\pm$$
which implies that
$$\pm i \iota_X\iota_Y d\omega_\pm + \iota_X\iota_Y \gamma = 0 \iff \mp i \left(d\omega_\pm\right)^{2,1}_\pm = \mp i \partial_\pm \omega_\pm= \gamma^{(2,1) + (3,0)}_\pm$$
Since $\omega_\pm$ and $\gamma$ are both real forms, we have
$$\gamma = \mp i (\partial_\pm \omega_\pm - \bar\partial_\pm \omega_\pm) = \pm d^c_\pm \omega_\pm$$
This recovers the correspondence (\cite{Gualtieri04, Gualtieri10}) of generalized K\"ahler structures with bi-Hermitian structures with torsion.

The Dorfman bracket induces the Bismut connections (Hitchin \cite{Hitchin06}). Let $X, Y \in TM$, then
$$\nabla^\pm_X Y := a\left\{\left[(X\mp g(X))*(Y\pm g(Y))\right]^{\pm}\right\} = \nabla_X Y \pm \frac{1}{2}g^{-1}\iota_X\iota_Y \gamma$$
where $\bullet^\pm$ denotes projection to $C_\pm$ respectively.
It follows that $\nabla^\pm$ preserve $g$ and their torsions are given by $\pm g^{-1} \iota_X\iota_Y \gamma$.
Since $L_\pm$ are involutive with respect to the Dorfman bracket, we see that $\ell_\pm$ as well as their conjugates are also involutive. In particular
$$Y \in T^+_{1,0} M \Longrightarrow \nabla^+_X Y \in a \left\{((\ell_- \dsum \bar \ell_-) * \ell_+)^+\right\} \subseteq a\{(\ell_- \dsum \ell_+ \dsum \bar \ell_-)^+\} = a(\ell_+) = T^+_{1,0} M$$
It follows that $\nabla^+$ preserves $I_+$. Similarly, $\nabla^-$ preserves $I_-$. 
We will use $\nabla^\pm$ to denote the induced Bismut connection on $T^*M$:
$$X \alpha(Y) = (\nabla^\pm_X \alpha)(Y) + \alpha(\nabla^\pm_X Y) \Longrightarrow \nabla^\pm \alpha = \nabla\alpha \pm \frac{1}{2}\iota_{\alpha^\flat} \gamma$$
Since $\nabla^\pm$ preserve $I_\pm$ respectively, they preserve $T^{0,1}_\pm M$ respectively as well. 

\subsection{Holomorphic reduction}\label{subsect:holomorphicreduction}
Consider the decomposition of the Courant algebroid $\TT_\CC M$ defined by a generalized K\"ahler structure:
\begin{equation}\label{eq:genkahlerdecomposition}
 \TT_\CC M = \ell_+ \dsum \ell_- \dsum \bar\ell_+ \dsum \bar \ell_-
\end{equation}
Reduction with respect to $\bar\ell_\pm$ defines $I_\pm$-holomorphic Courant algebroids
$$0 \to T^{1,0}_\pm M \to \mc E_\pm \to T_{1,0}^\pm M \to 0$$
For example, the $I_+$-holomorphic vector bundle $\mc E_+$ is given by
$$\mc E_+ := \frac{(\bar\ell_+)^\perp}{\bar\ell_+} = \frac{\bar\ell_+ \dsum \bar\ell_- \dsum \ell_-}{\bar\ell_+} \cong \bar\ell_- \dsum \ell_-$$
where the $I_+$-holomorphic structure is induced by the Dorfman bracket
\begin{equation}\label{eq:inducedholostructure}
 \bar\partial_{+,X} \mf Y_A = (\mf X * \mf Y)_A\text{ for } \mf X \in \bar\ell_+, X = a(\mf X) \in T_{0,1}^+ M, \mf Y \in (\bar\ell_+)^\perp
\end{equation}
and $\bullet_A$ denotes the equivalence class.
Under this reduction, the Lie algebroid $\bar L_+$ induces an $I_+$-holomorphic Lie subalgebroid in $\mc E_+$:
$$\bar L_+ = \bar\ell_+ \dsum \bar\ell_- \Longrightarrow \mc A_+ := \frac{\bar\ell_+ \dsum \bar\ell_-}{\bar\ell_+} \cong \bar \ell_- \subset \mc E_+$$
Similarly, reduction by $\bar\ell_-$ defines the $I_-$-holomorphic vector bundle $\mc E_- \cong \bar\ell_+ \dsum \ell_+$ and $\bar L_+$ induces an $I_-$-holomorphic Lie subalgebroid $\mc A_- \cong \bar\ell_+$ in $\mc E_-$. We caution here that the isomorphisms in this paragraph are only isomorphisms of complex vector bundles.

We consider the $I_+$-holomorphic Lie algebroid $\mc A_+$. For  $\mf X_\pm \in C^\infty(\bar\ell_\pm)$, write $\mf X = \mf X_+ + \mf X_-$.
By definition, we have $\mf X_A = (\mf X_-)_A$. The bracket \emph{does not} descent to $C^\infty(\mc A_+)$, because
$$[\mf X_A, \mf Y_A] = \left(\mf X * \mf Y\right)_A = \left(\mf X_- * \mf Y_- + \mf X_+ * \mf Y_- + \mf X_- * \mf Y_+ \right)_A$$
becomes contradictory, for example, when $\mf X_- = \mf Y_+ = 0$. On the other hand, when restricted to the $I_+$-holomorphic sections, we have
$$[\mf X_A, \mf Y_A] = (\mf X_- * \mf Y_-)_A = ([\mf X_-, \mf Y_-])_A$$
Thus, $\bar\ell_\mp$ are the \emph{smooth} Lie algebroids underlying $\mc A_\pm$.

\subsection{Lie algebroid modules} \label{subsect:Liealgebroidcohomology}
Let $L$ be a smooth Lie algebroid. The \emph{Lie algebroid de Rham complex} $(\Omega^\bullet(L), d_{L})$ is 
$$0 \to C^\infty(M) \xto{d_{L}} C^\infty(L^*) \xto{d_{L}} C^\infty(\wedge^2 L^*) \xto{d_{L}} \ldots$$
where it is customary to denote $\Omega^k(L) = C^\infty(\wedge^k L^*)$. 
The differential $d_{L}$ is defined in the standard fashion. For $\alpha \in \Omega^k(L)$, and $s_0, \ldots, s_k \in C^\infty(L)$ with $X_i = \a(s_i)$,
\begin{equation} \label{eq:Liealgebroiddiffdefin}
 \begin{split}
  d_{L} \alpha (s_0,\ldots, s_k) := & \sum_{j = 0}^k (-1)^j X_j \alpha(s_0, \ldots, \hat s_j, \ldots, s_k) \\
  & + \sum_{i < j} (-1)^{i+j} \alpha(s_i * s_j, s_0, \ldots, \hat s_i, \ldots, \hat s_j, \ldots, s_k)
 \end{split}
\end{equation}
The homology of this complex is the \emph{Lie algebroid cohomology}, denoted $H^k(L)$. 

Let $V$ be a vector bundle, an \emph{$L$-connection} is a derivation $D : V \to V \tensor L^*$:
$$D(fv) = d_Lf \tensor v + f Dv \text{ for } f \in C^\infty(M)$$
It is naturally extended to $\Omega^k(L; V) := C^\infty(V \tensor \wedge^k L^*)$ using $d_{L}$:
\begin{equation} \label{eq:Liealgebroidconnectionextension}
 D(\alpha \tensor v) = d_L \alpha \tensor v + (-1)^k \alpha \tensor Dv
\end{equation}
We say that $(V, D)$ (or simply $V$) is an \emph{$L$-module} if $D\circ D = 0$, i.e. $D$ is flat. In this case, there is a induced complex $(\Omega^*(V; L), D)$, whose homology is denoted $H^*(V; L)$.

For a holomorphic Lie algebroid $\mc L$, we have correspondingly the \emph{holomorphic Lie algebroid de Rham complex} $(\wedge^* \mc L, \partial_\mc L)$. For local holomorphic sections $\alpha \in \wedge^k \mc L^*$ and $s_0, \ldots, s_k \in \mc L$, with $X_i = \a(s_i)$ local holomorphic vector fields, the righthand side of \eqref{eq:Liealgebroiddiffdefin} defines $\partial_\mc L \alpha(s_0, \ldots, s_k)$. The hyperhomology is 
denoted $\HH^*(\mc L)$. Similarly, we can define the notion of \emph{$\mc L$-connection} on a holomorphic vector bundle $V$, as well as the notion of \emph{$\mc L$-modules}\footnote{A more general notion is a \emph{sheaf of $\mc L$-modules} (cf. Tortella \cite{Tortella11}). The $\mc L$-modules considered here are the locally free ones.}. For an $\mc L$-module $V$, the hyperhomology of the induced complex $(V \tensor \wedge^* \mc L^*, \partial_\mc L)$ is denoted $\HH^*(V;\mc L)$.

\begin{remark}\label{remark:Atiyahclass}
 \rm{
  For a holomorphic vector bundle $V$, the classical Atiyah class $\alpha(V) \in H^1(T_{1,0}^*M \tensor \End(V))$ is the obstruction of existence of a holomorphic connection (Atiyah \cite{Atiyah57}).
  One can show that the existence of an $\mc L$-connection is obstructed by the corresponding Atiyah class $a^*\alpha(V) \in H^1(\mc L^* \tensor \End(V))$. When $V$ is a line bundle, $a^*\alpha(V) \in H^1(\mc L^*)$. For an $\mc L$-module $V$, we have
  $$a^*\alpha(V) = 0$$
  We note that in general, this does not imply that $\alpha(V) = 0$.
 }
\end{remark}

\begin{defn}\label{defn:LiealgebroidPicardgroup}
 Let $\mc L$ be a holomorphic Lie algebroid. Two $\mc L$-modules $V_1$ and $V_2$ are \emph{equivalent} if there exists holomorphic bundle isomorphism $f: V_1 \to V_2$ covering identity that preserves the $\mc L$-connections. The \emph{$\mc L$-Picard group} $\Pic^{\mc L}(M)$ consists of the equivalence classes of $\mc L$-modules of rank $1$, where the group structure is given by tensor product.
\end{defn}

The group $\Pic^{\mc L}(M)$ is obviously abelian.
\begin{lemma}\label{lemma:Picardhomomorphism}
 The map $\Pic^{\mc L}(M) \to \Pic(M)$ forgetting the $\mc L$-connection is a group homomorphism, whose kernel $\Pic^{\mc L}_0(M)$ is the subgroup where the underlying holomorphic line bundle is trivial.
 \qed
\end{lemma}

\subsection{$\JJ_+$-holomorphic vector bundles}\label{subsect:Jplusholomorphicbundles}
Let $\JJ$ be a generalized complex structure, with the $i$-eigensubbundle $L$.
A \emph{$\JJ$-holomorphic vector bundle} $V \to M$ is an $\bar L$-module, i.e. a complex vector bundle with a flat $\bar L$-connection:
$$D : C^\infty(V) \to C^\infty(V \tensor \bar L^*), D(fv) = d_{\bar L} f \tensor v + f Dv$$
such that $D\circ D = 0$. 
The homology of the associated complex $(\Omega^*(V; \bar L), D)$ is denoted $H^*(V; \JJ)$. We say that a section $v \in \Omega^0(V)$ is \emph{$\JJ$-holomorphic} if $Dv = 0$. Similar to the classical case, $H^0(V; \JJ)$ consists of global holomorphic sections of $V$.

Suppose now that $\JJ = \JJ_+$ is part of a generalized K\"ahler structure.
Since $\bar L = \bar \ell_+ \dsum \bar \ell_-$, $D$ naturally decomposes as $D = D_+ + D_-$:
$$D_\pm : C^\infty(V) \to C^\infty(V \tensor \bar \ell_\pm^*), D_\pm(fv) = d_{\bar \ell_\pm} f \tensor v + f D_\pm v$$
which implies that $D_\pm$ are $\bar \ell_\pm$-connections on $V$. We use the same notations $D_\pm$ to denote their extensions to $\Omega^*(V; \bar \ell_\pm)$ by \eqref{eq:Liealgebroidconnectionextension}.
\begin{prop}\label{prop:genholoisbiholo}
 Let $\JJ = \JJ_+$ be part of a generalized K\"ahler structure. A $\JJ$-holomorphic vector bundle is an $I_\pm$-holomorphic vector bundle  on $(M, I_+, I_-)$. A section $v \in \Omega^0(V)$ is $\JJ$-holomorphic iff it is a holomorphic section in both $I_\pm$-holomorphic structures.
\end{prop}
{\it Proof:} 
 Let $s_\pm \in \bar\ell_\pm$ respectively and write $s = s_+ + s_- \in \bar L$, then for $v \in V$, we have
  $$D_{\pm, s_\pm} (v) =  D_{s_\pm}(v)$$
 The flatness of $ D$ is equivalent to the following identity for any $s, t \in \bar L$ and $v \in V$:
 \begin{equation}\label{eq:flatnessdiffequation}
   D_s D_t(v) -  D_t D_s(v) -  D_{s*t}(v) = 0 
 \end{equation}
 By setting $s_\bullet = t_\bullet = 0$, $\bullet = \pm$, in the above equality we see that $ D_\pm$ are flat as well. Since $\bar \ell_\pm \cong T_{0,1, \pm} M$ as Lie algebroids, we see that $D_\pm$ defines $I_\pm$-holomorphic structures on $V$.
 The last statement is straightforward.
\qed

In \cite{Gualtieri10}, using the notion of Morita equivalence, it is shown that the category of $\JJ_+$-holomorphic vector bundles is equivalent to the category of $\mc A_\bullet$-modules, where $\bullet = +$ or $-$. We will show that this equivalence is in fact an isomorphism of categories (Proposition \ref{prop:isomorphiccategories}).

\subsection{Canonical line bundle} \label{subsect:spinorbundle}
Let $\JJ$ be a generalized complex structure.
The \emph{$\JJ$-canonical bundle} $U$ is the pure spinor line bundle defining the Dirac structure $\bar L$, i.e.
$$\bar L = \{\mf X | \mf X \cdot C^\infty(U) = 0\}$$
Let $\chi \in C^\infty(U)$ and $d_\gamma := d + \gamma \wedge$ the \emph{twisted de Rham differential}, it's shown in \cite{Gualtieri04} that the integrability of $\JJ$ is equivalent to $d_\gamma \chi \in C^\infty(L) \cdot C^\infty(U)$. Since the spinor action defines an isomorphism of $L \tensor U$ onto the image, $d_\gamma$ defines an $\bar L$-connection under the identification $L \cong \bar L^*$ by $2\<,\>$:
\begin{equation}\label{eq:canonicalbundleconnection}
 D_{U} : C^\infty(U) \to C^\infty(U \tensor \bar L^*) : D_{U, s} \chi := s \cdot d_\gamma \chi \text{ for } s \in C^\infty(\bar L)
\end{equation}
It turns out that $D_{U}$ is flat (\cite{Gualtieri04}) and $U$ is naturally a $\JJ$-holomorphic line bundle.

When $\JJ = \JJ_+$ is part of a generalized K\"ahler structure, the canonical line bundle generates the \emph{generalized Hodge decomposition} of the twisted de Rham complex $(\Omega^*(M), d_\gamma)$ via the spinor actions of $\JJ_\pm$ (\cite{Gualtieri0409} and Baraglia \cite{Baraglia12}). More precisely, let $U_{-n, 0} := C^\infty(U) \subset \Omega^*(M)$, then
$$U_{r,s} := \left(C^\infty(\wedge^p \ell_+) \tensor C^\infty(\wedge^q \ell_-)\right) \cdot U_{-n,0} \subset \Omega^*(M)$$
with $r = p+q - n$ and $s = p-q$, by the spinor action of $C^\infty(\ell_\pm)$ on $U_{-n,0}$. The differential $d_\gamma$ decomposes as 
\begin{equation}\label{eq:dgammadecomposition}
 d_\gamma = \partial_+^\gamma + \bar\partial_+^\gamma \text{ with } \partial_+^\gamma := \delta_+^\gamma + \delta_-^\gamma \text{ and } \bar\partial_+^\gamma := \bar\delta_+^\gamma + \bar\delta_-^\gamma
\end{equation}
where $\bar\delta_\pm^\gamma : U_{r,s} \to U_{r+1, s\pm 1}$ and $\delta_\pm^\gamma : U_{r,s} \to U_{r-1, s\mp 1}$. In terms of the gradings given by $p$ and $q$, we have
\begin{equation}\label{eq:twisteddiffdecompositions}
 \begin{split}
  \bar\delta_+^\gamma & : \left(C^\infty(\wedge^p \ell_+) \tensor C^\infty(\wedge^q \ell_-)\right) \cdot U_{-n,0} \to \left(C^\infty(\wedge^{p+1} \ell_+) \tensor C^\infty(\wedge^q \ell_-)\right) \cdot U_{-n,0} \\
  \bar\delta_-^\gamma & : \left(C^\infty(\wedge^p \ell_+) \tensor C^\infty(\wedge^q \ell_-)\right) \cdot U_{-n,0} \to \left(C^\infty(\wedge^p \ell_+) \tensor C^\infty(\wedge^{q+1} \ell_-)\right) \cdot U_{-n,0} \\
  \delta_+^\gamma & : \left(C^\infty(\wedge^p \ell_+) \tensor C^\infty(\wedge^q \ell_-)\right) \cdot U_{-n,0} \to \left(C^\infty(\wedge^{p-1} \ell_+) \tensor C^\infty(\wedge^q \ell_-)\right) \cdot U_{-n,0} \\
  \delta_-^\gamma & : \left(C^\infty(\wedge^p \ell_+) \tensor C^\infty(\wedge^q \ell_-)\right) \cdot U_{-n,0} \to \left(C^\infty(\wedge^p \ell_+) \tensor C^\infty(\wedge^{q-1} \ell_-)\right) \cdot U_{-n,0} \\
 \end{split}
\end{equation}
In \cite{Gualtieri0409}, it's shown that the Laplacians of the various differentials coincide up to scalars, which gives the generalized Hodge decomposition of the twisted de Rham cohomology $H_\gamma^*(M)$.

\section{A spectral sequence} \label{sect:liealgebroidspectralsequence}
We describe in this section a computational tool, the \emph{BH-to-GK} spectral sequence, relating the cohomology groups in the bi-Hermitian data to those in the generalized K\"ahler data. 
\subsection{The double complex} \label{subsect:doublecomplex}

Let's consider the decomposition $L := L_+ = \ell_+ \dsum \ell_-$. This is \emph{not} a direct sum as Lie algebroids. Although the brackets on $L$ as well as the summands are given by the restriction of the Dorfman bracket on $\TT M$, in general
$$\ell_+ * \ell_- \neq 0$$
Analysing this non-vanishing leads to the spectral sequence.

We have the bundle isomorphism 
$$\wedge^*\bar L^* \cong \Dsum_{p,q} \wedge^p\bar \ell_+^* \tensor \wedge^q \bar\ell_-^* \cong \Dsum_{p,q} \wedge^p T_+^{0,1} M \tensor \wedge^q T^{0,1}_- M$$
which induces the isomorphism as linear spaces
$$\Omega^*(\bar L) \cong \Dsum_{p,q} \Omega^{p,q}(\bar L) \text{ where } \Omega^{p,q}(\bar L) = \Omega^p(\bar\ell_+) \tensor \Omega^q(\bar\ell_-) \cong \Omega_+^{0,p}(M) \tensor \Omega_-^{0,q}(M)$$
Let $\alpha \in \Omega_+^{0,p}(M)$ and $\beta \in \Omega_-^{0,q}(M)$. The same notations are used for the corresponding elements in $\Omega^*(\bar\ell_\pm)$. Then $\alpha \wedge \beta \in \Omega^{p,q}(\bar L)$. By definition, we have
$$d_{\bar L}(\alpha \wedge \beta) = d_{\bar L} \alpha \wedge \beta + (-1)^p \alpha \wedge d_{\bar L} \beta$$
\begin{lemma}\label{lemma:compareLiealgebroiddiff}
 Let $\bar\partial^\mp : \Omega^*(M) \xto{\nabla^\pm} \Omega^*(M) \tensor \Omega^1(M) \xto{} \Omega^*(M) \tensor \Omega_\mp^{0,1}(M)$. Then 
 $$\left.d_{\bar L}\right|_{\Omega^{p,0}(\bar L)} = \bar\partial_+ + \bar\partial^- \text{ while } \left.d_{\bar L}\right|_{\Omega^{0,q}(\bar L)} = \bar\partial_- + \bar\partial^+$$
\end{lemma}
{\it Proof:}
 We prove the $+$-case, by directly comparing the two differentials $d_{\bar L}$ and $d_{\bar\ell_+}$. Consider the sections $s_i \in \bar\ell_+$ for $0 < i \varleq p$, $r_\pm \in \bar\ell_\pm$ and $r = r_+ + r_- \in L$. Suppose that $X_i = a(s_i)$, $Z_\pm = a(r_\pm)$ and $Z = a(r)$. We compute for $\alpha \in \Omega_+^{0,p}(M)$
 \begin{equation*}
  \begin{split}
   & \left(d_{\bar L}\alpha\right)(r, s_1, \ldots, s_p) = Z \alpha(s_1, \ldots, s_p) + \sum_{i = 1}^p (-1)^i X_i\alpha(r_+, \ldots, \hat s_i, \ldots, s_p)\\
   & + \sum_{i = 1}^p (-1)^i\alpha((r*s_i)^+, \ldots, \hat s_i, \ldots, s_p) + \sum_{i<j}(-1)^{i+j}\alpha(s_i *s_j, r_+, \ldots, \hat s_i, \ldots, \hat s_j, \ldots, s_p) \\
   = & \left(d_{\bar \ell_+}\alpha\right)(r_+, s_1, \ldots, s_p) + Z_-\alpha(s_1, \ldots, s_p) + \sum_{i = 1}^p (-1)^{i} \alpha((r_-*s_i)^+, \ldots, \hat s_i, \ldots, s_p)
  \end{split}
 \end{equation*}
 When we apply the identification of $\ell_\pm$ with the holomorphic tangent bundles with respect to $I_\pm$, we see that the first term corresponds to $\bar\partial_+ \alpha$ and last two terms give precisely $\bar\partial^-_{r_-} \alpha$.
\qed

We are ready to describe the BH-to-GK spectral sequence. By Lemma \ref{lemma:compareLiealgebroiddiff}, on $\Omega^{p,q}(\bar L)$, we define
\begin{equation}\label{eq:deltaoperators}
 \bar\delta_+ : = \bar\partial_+ + \bar\partial^+ : \Omega^{p,q}(\bar L) \to \Omega^{p+1, q}(\bar L) \text{ and } \bar\delta_- : = \bar\partial_- + \bar\partial^- : \Omega^{p,q}(\bar L) \to \Omega^{p, q+1}(\bar L)
\end{equation}
Then $d_{\bar L} = \bar\delta_+ + \bar\delta_-$. 
For example, with $\alpha \in \Omega^{p,0}(\bar L)$ and $\beta \in \Omega^{0,q}(\bar L)$,
$$\bar\delta_+(\alpha \wedge \beta) = \bar\partial_+\alpha \wedge \beta + (-1)^p \alpha \wedge \bar\partial^+ \beta$$
In particular, for $\alpha \in \Omega^{p,0}(\bar L) \cong \Omega_+^{0,p}(M)$, we have
$$\bar\delta_+ \alpha = \bar\partial_+ \alpha \text{ and } \bar\delta_- \alpha = \bar\partial^- \alpha$$
and similarly for $\beta \in \Omega^{0,q}(\bar L) \cong \Omega_-^{0,q}(M)$, we have
$$\bar\delta_+ \beta = \bar\partial^+ \beta \text{ and } \bar\delta_- \beta = \bar\partial_- \beta$$
\begin{prop}\label{prop:BHtoGK}
 The spectral sequences associated with the double complex $(\Omega^{p,q}(\bar L), \bar\delta_+, \bar\delta_-)$ compute the Lie algebroid cohomology $H^*(\bar L)$.
\end{prop}
{\it Proof:}
 From $d_{\bar L}^2 = 0$, we get $\bar\delta_+^2 = 0$, $\left[\bar\delta_+, \bar\delta_-\right] = 0$ and $\bar\delta_-^2 = 0$. Thus, $(\Omega^{p,q}(\bar L), \bar\delta_+, \bar\delta_-)$ is indeed a double complex. By definition, the corresponding spectral sequences compute $H^*(\bar L)$.
\qed

\begin{remark}\label{remark:GualtieriGKpaper}
 \rm{
  The double complex $(\Omega^{p,q}(\bar L), \bar\delta_+, \bar\delta_-)$ appeared already in the proof of Proposition 3.11 in \cite{Gualtieri10}.
  Here we identify explicitly the differentials using the bi-Hermitian data. 
 }
\end{remark}

\begin{defn}\label{defn:BHtoGK}
 The spectral sequences associated to the double complex $(\Omega^{p,q}(\bar L), \bar \delta_+, \bar\delta_-)$ are the \emph{BH-to-GK spectral sequences} for the generalized K\"ahler manifold $(M, \gamma; \JJ_\pm, \sigma)$.
\end{defn}

The \emph{BH-to-GK} spectral sequences are related to the \emph{holomorphic reduction} (\cite{Gualtieri10}, and \S \ref{subsect:holomorphicreduction}).

\begin{lemma}\label{lemma:explicitreducedholomorphicstructure}
 The operators $\bar\delta_\pm$ on $\bar\ell_\mp^*$ defined in \eqref{eq:deltaoperators} is induced by the operators $\bar\partial_{\pm}$ on $\mc A_\pm$ defined in \eqref{eq:inducedholostructure} under the natural isomorphism $\bar\ell_\mp \cong \mc A_\pm$.
\end{lemma}
{\it Proof:}
 This follows essentially by comparing the definitions in \eqref{eq:deltaoperators} and \eqref{eq:inducedholostructure}. We write down the proof for $\mc A_+$. Let $\alpha_A \in \mc A_+^*$, with $\alpha$ the corresponding element in $\Omega^{0,1}(\bar L) = \bar\ell_-^*$ under the natural isomorphism. Let $\mf X_+ \in \bar\ell_+$ and $\mf Y_- \in \bar\ell_-$, with $X = a(X_+) \in T_{0,1}^+ M$ and $(\mf Y_-)_A \in \mc A_+$. We compute
 \begin{equation*}
  \begin{split}
   & \left(\bar\partial_{+, X} \alpha_A\right) ((\mf Y_-)_A) = X \alpha_A((\mf Y_-)_A) - \alpha_A\left(\bar\delta'_{+,X} (\mf Y_-)_A\right) \\
   \mapsto\hspace{0.1in} & X\alpha(\mf Y_-) - \alpha((\mf X_+ * \mf Y_-)^-) = \left(\bar\delta_{+,X}\alpha\right)(\mf Y_-)
  \end{split}
 \end{equation*}
 Thus $\bar\partial_+$ and $\bar \delta_+$ coincide on $\mc A_+^* \cong \bar\ell_-^*$.
\qed

We can now say a little more about the spectral sequences (cf. \cite{Gualtieri10}). Look vertically, Lemma \ref{lemma:compareLiealgebroiddiff} implies that the double complex $(\Omega^{p,q}(\bar L), \bar\delta_+, \bar\delta_-)$ gives a resolution of the holomorphic de Rham complex of $\mc A_-$. The first page of the \emph{first} spectral sequence $({}_IE^{p,q}_r, d_r)$ consists of $\bar\partial_-$-Dolbeault cohomologies:
$${}_IE^{p,q}_1 = H^{0,q}_{\bar\partial_-}(\wedge^p \mc A_-^*) \text{ with } d_1 \text{ induced by } \bar\delta_+$$ 
Similarly, the first page of the the \emph{second} spectral sequence is
$${}_{II}E^{p,q}_1 = H^{0,q}_{\bar\partial_+}(\wedge^p \mc A_+^*) \text{ with } d_1' \text{ induced by } \bar\delta_-$$

\begin{prop}\label{prop:isomorphiccohomology}
 $\HH^*(\mc A_+) \cong \HH^*(\mc A_-) \cong H^*(\bar L_+)$.
 \qed
\end{prop}

\begin{remark}\label{remark:matchedpair}
 \rm{
  The holomorphic Lie algebroids $\mc A_\pm$ form a \emph{matched pair} (cf. \cite{Gualtieri10}, Laurent-Gengoux-Sti\'enon-Xu \cite{LaurentGengouxStienonXu08}, Mackenzie \cite{Mackenzie07}, Mokri \cite{Mokri97} and Lu \cite{Lu97}). Proposition \ref{prop:isomorphiccohomology} corresponds to, for example, Theorem $4.11$ of \cite{LaurentGengouxStienonXu08}.
 }
\end{remark}

\begin{example}\label{exple:Kahlerspectralsequence}
 \rm{
  Let $(M, J, \omega)$ be a K\"ahler manifold, which can be seen as a generalized K\"ahler manifold as following:
  $$\JJ_+ = \begin{pmatrix} J & 0 \\ 0 & -J^*\end{pmatrix} \text{ and } \JJ_- = \begin{pmatrix}0 & -\omega^{-1} \\ \omega & 0\end{pmatrix}$$
  It follows that $\bar L_+ = T_{0,1}M \dsum T^{1,0} M$, and the decomposition into $\bar \ell_\pm$ is given by the metric $g$:
  $$\bar \ell_\pm = \{X \pm g(X) : X \in T_{0,1} M\}$$
  Thus, $I_+ = I_- = J$. The holomorphic reduction in this case endows, for example on $\bar \ell_-$, the holomorphic structure given by
  $$\bar\partial_X (Y - g(Y)) \equiv (X+ g(X)) * (Y - g(Y))  \equiv \nabla_X Y - g(\nabla_X Y) \mod \bar \ell_+ \text{ for } X, Y \in C^\infty(T_{0,1}M)$$
  where $\nabla$ is the Levi-Civita connection. By the K\"ahler condition, $\nabla$ preserves $J$, which implies that $\nabla_X Y \in C^\infty(T_{0,1} M)$ as well. Let $\nabla^{0,1}$ denote the $(0,1)$-component of $\nabla$, then we showed that
  $$\bar\partial = \nabla^{0,1} \text{ on } C^\infty(\bar \ell_-) \cong C^\infty(T_{0,1} M)$$
  Lemma \ref{lemma:explicitreducedholomorphicstructure} implies that $\nabla^{0,1}$ induces the operator $\bar\delta_+$ in the corresponding double complex $\Omega^{p,q}(\bar L)$. The metric $g$ defines an isomorphism of vector bundles:
  $$\bar \ell_- \to T_{0,1} M \xto{g} T^{1,0} M \Rightarrow \bar \ell_-^* \cong T_{1,0} M$$
  The induced holomorphic structure on $T_{1,0} M$ coincides with the standard holomorphic structure as the holomorphic tangent bundle. By Proposition \ref{prop:isomorphiccohomology}, we get $H^*(\bar L_+) \cong \HH^*(T^{1,0} M)$. In particular, the first page of the spectral sequence is given by the Dolbeault cohomology of the sheaves of the holomorphic multivectors:
  $$E_1^{p,q} = H^{0,q}_{\bar \partial}(\wedge^p T_{1,0} M)$$
  Since in this case, $I_+ = I_-$, $\ell_+$ can be used in the above arguments and arrive at the same conclusion.
 }
\end{example}

\subsection{$\JJ_+$-holomorphic vector bundles} \label{subsect:genholobundlespectralsequence}
The \emph{BH-to-GK} spectral sequence induces a corresponding spectral sequence for $\JJ_+$-holomorphic vector bundles. Again, we write $L := L_+$. 

Let $(V, D)$ be a $\JJ$-holomorphic vector bundle. Then we have the decomposition $D = D_+ + D_-$ (see \S\ref{subsect:Jplusholomorphicbundles}) into $\bar \ell_\pm$-connections. The extension of the decomposition to $\Omega^*(V; \bar L)$ is induced by the corresponding decompositions $\displaystyle{\Omega^*(\bar L) = \dsum_{p,q} \Omega^{p,q}(\bar L)}$ and $d_{\bar L} = \bar\delta_+ + \bar\delta_-$:
$$\Omega^* (V; \bar L) = \Dsum_{p,q} \Omega^{p,q}(V; \bar L) \text{ where } \Omega^{p,q}(V; \bar L) = C^\infty(V \tensor \wedge^p \bar\ell_+^* \tensor \wedge^q \bar \ell_-^*)$$
We write down the extension of $D_+$:
\begin{equation}\label{eq:extendingdifferential}
 D_+ : \Omega^{p,q}(V; \bar L) \to \Omega^{p+1, q}(V, \bar L) : D_+(\alpha \tensor v) = \bar\delta_+\alpha \tensor v + (-1)^{p+q}\alpha \tensor D_+ v
\end{equation}
Since $D^2 = 0$, it is straightforward to verify that $(\Omega^{p,q}(V; \bar L), D_+, D_-)$ is a double complex. The equation $[D_+, D_-] = 0$ is called the \emph{commutation relation}. 
\begin{defn}\label{defn:BHtoGKbundles}
 The spectral sequences associated to the double complex $(\Omega^{p,q}(V; \bar L), D_+, D_-)$ is the \emph{BH-to-GK spectral sequences} for the $\JJ_+$-holomorphic vector bundle $V$.
\end{defn}

The commutation relation helps us understand the equivalence of categories in \cite{Gualtieri10}. We start with the following improvement on Proposition \ref{prop:genholoisbiholo}.
\begin{lemma}\label{lemma:commutationrelationdefinesgenholo}
 A complex vector bundle $V$ is $\JJ_+$-holomorphic iff it is an $I_\pm$-holomorphic vector bundle and the commutation relation holds for the induced differentials $D_\pm$.
\end{lemma}
{\it Proof:}
 We show the $\Longleftarrow$ direction. Let $\bar\partial_\pm$ denote the $I_\pm$-holomorphic structures on $V$, i.e.
 $$\bar\partial_\pm : C^\infty(V) \to C^\infty(V \tensor T^{0,1}_\pm M) \text{ and } \bar\partial_\pm \circ \bar\partial_\pm = 0$$
 Composing with the isomorphisms $T^{0,1}_\pm M \cong \bar \ell_\pm^*$ gives the operators $D_\pm : C^\infty(V) \to C^\infty(V \tensor \bar\ell_\pm^*)$. Let $D = D_+ + D_-$. Since $d_{\bar L} = d_{\bar\ell_+} + d_{\bar\ell_-}$ on functions, we see that $D$ is an $\bar L$-connection. Extend $D_\pm$ to $\Omega^{p,q}(V;\bar L)$ by \eqref{eq:extendingdifferential}, then $D = D_+ + D_-$ is extended to $\Omega^*(V; \bar L)$. We compute
 $$D\circ D = (D_+ + D_-) \circ (D_+ + D_-) = D_+ \circ D_+ + [D_+, D_-] + D_- \circ D_- = 0$$
 which implies the that $V$ is a $\JJ_+$-holomorphic bundle.
\qed

\begin{lemma}\label{lemma:holotomodule}
 Let $(V, D)$ be a $\JJ_+$-holomorphic vector bundle. Then $D_\mp$ define structures of $\mc A_\pm$-modules on the induced $I_\pm$-holomorphic vector bundles $(V, D_\pm)$.
\end{lemma}
{\it Proof:}
 We prove the case of $V_+$. The $I_+$-holomorphic structure on $V$ induces one on $V \tensor \mc A_+^* \cong V \tensor \bar\ell_-^*$, both of which coincide with $D_+$ (cf. Lemma \ref{lemma:compareLiealgebroiddiff} and \eqref{eq:extendingdifferential}):
 $$D_+ : C^\infty(V) \to C^\infty(V \tensor \bar\ell_+^*) \Longrightarrow D_+ : C^\infty(V \tensor \bar\ell_-^*) \to C^\infty(V \tensor \bar\ell_-^* \tensor \bar\ell_+^*)$$
 whose (local) kernels consist of (local) $I_+$-holomorphic sections.
 For $D_- : C^\infty(V) \to C^\infty(V \tensor \bar\ell_-^*)$, the commutation relation $[D_+, D_-] = 0$ implies that when $v$ is a local $I_+$-holomorphic section of $V$, $D_-(v)$ is a local $I_+$-holomorphic section of $V \tensor \bar\ell^*_-$. Thus, we obtain $\partial_{\mc A_+} : V \to V \tensor \mc A_+^*$. 
\qed

Similar to Proposition \ref{prop:isomorphiccohomology}, we have
\begin{prop}\label{prop:isomorphicbundlecohomology}
 $\HH^*(V; \mc A_+) \cong \HH^*(V; \mc A_-) \cong H^*(V; \JJ_+)$. Moreover, the first page of the \emph{first} spectral sequence $({}_IE^{p,q}_r, d_r)$ is
  $${}_IE^{p,q}_1 = H^{0,q}_{\bar\partial_-}(V \tensor \wedge^p \mc A_-^*) \text{ with } d_1 \text{ induced by } D_+$$ 
 Similarly, the first page of the the \emph{second} spectral sequence is
  $${}_{II}E^{p,q}_1 = H^{0,q}_{\bar\partial_+}(V \tensor \wedge^p \mc A_+^*) \text{ with } d_1' \text{ induced by } D_-$$
\end{prop}
{\it Proof:}
 The vertical and horizontal sequences in the double complex $(\Omega^{p,q}(V; \bar L), D_+, D_-)$ are respectively the resolutions of the holomorphic de Rham complexes of $\mc A_\mp$
 (cf. \cite{Wells79}). 
\qed

\begin{remark}\label{remark:refertoXupaper}
 \rm{
  As in Remark \ref{remark:matchedpair}, Proposition \ref{prop:isomorphicbundlecohomology} corresponds to Theorem $4.19$ of \cite{LaurentGengouxStienonXu08}.
 }
\end{remark}

Suppose now that $V$ is an $\mc A_+$-module (the case of $\mc A_-$-module is similar). Namely
\begin{enumerate}
 \item $V$ is $I_+$-holomorphic, with the $I_+$-holomorphic structure $\bar\partial_+ : C^\infty(V) \to C^\infty(V \tensor T^{0,1}_+ M)$
 \item together with a flat $\mc A_+$-connection, i.e. $\partial_{\mc A_+} : V \to V \tensor \mc A_+^*$ and $\partial_{\mc A_+} \circ \partial_{\mc A_+} = 0$.
\end{enumerate}
We recover the structure of a $\JJ_+$-holomorphic vector bundle on $V$. Since $\bar\ell_-$ is the underlying smooth Lie algebroid of $\mc A_+$ and $\bar\ell_- \cong T_{0,1}^- M$ as smooth Lie algebroids, $\partial_{\mc A_+}$ induces an $I_-$-holomorphic structure $\bar\partial_-$ on $V$:
\begin{enumerate}
 \item For a local holomorphic section $v$, $\bar\partial_- v := \partial_{\mc A_+} v$ under the isomorphism $\mc A_+ \cong T_{0,1}^- M$
 \item For a local smooth function $f$, $\bar\partial_- (fv) := \bar\partial_- f \tensor v + f\bar\partial_- v$
\end{enumerate}
This defines $\bar\partial_-$ for all local smooth sections of $V$. For the same $f$ and $v$ as above, we have 
$$\bar\partial_- \circ \bar\partial_- (fv) = f \partial_{\mc A_+} \circ \partial_{\mc A_+} v = 0$$
i.e., $V$ is an $I_\pm$-holomorphic bundle. By Lemma \ref{lemma:commutationrelationdefinesgenholo}, we need to verify the commutation relation.
\begin{lemma}\label{lemma:moduletoholo}
 Extend $\bar\partial_\pm$ to $D_\pm$ on $\Omega^{p,q}(V; \bar L)$ using \eqref{eq:extendingdifferential}, then the commutation relation holds.
\end{lemma}
{\it Proof:}
 Written in terms of covariant derivatives, the commutation relation becomes
  $$D_{+, s} D_{-,t} (v) - D_{-, \bar\delta_{+,s} t} (v) - D_{-, t} D_{+, s} (v) + D_{+, \bar\delta_{-,t} s} (v) = 0 \text{ for } s \in \bar\ell_+, t \in \bar\ell_-$$
 The equation $[\bar\delta_+, \bar\delta_-] = 0$ implies that the left hand side is tensorial in all variables. The situation can be simplified as following. 
 \begin{enumerate}
  \item We only need to consider when $v$ is a $I_+$-holomorphic section of $V$, i.e. $D_+ (v) = 0$.
  \item We may take $s$ to be $I_-$-holomorphic and $t$ to be $I_+$-holomorphic, i.e. 
   $$\bar\delta_{+,s} t = 0 \text{ and } \bar\delta_{-,t} s = 0$$
 \end{enumerate}
 The remaining term is $D_{+,s} D_{-,t}(v) = D_{+,s} \partial_{\mc A_+, t}(v) = 0$.
\qed

The two Lemmata \ref{lemma:holotomodule} and \ref{lemma:moduletoholo} imply that a complex vector bundle $V$ is $\JJ_+$-holomorphic iff it is an $\mc A_+$-module iff it is an $\mc A_-$-module. Going through the construction, it is fairly clear that they are functorial. In summary, we obtain a slight refinement of the statement in \cite{Gualtieri10} concerning equivalence of categories (Corollary $3.15$ \emph{loc. cit.}).
\begin{prop}\label{prop:isomorphiccategories}
 On a generalized K\"ahler manifold, the categories of $\JJ_+$-holomorphic vector bundles and (locally free) $\mc A_\pm$-modules are isomorphic to each other.
 \qed
\end{prop}

\begin{remark}\label{remark:recoverXu}
 \rm{
  The correspondence of $\JJ_+$-holomorphic vector bundles with $\mc A_\pm$-modules follows also from Lemma $4.16$ of \cite{LaurentGengouxStienonXu08}. Our proof is basically an adaptation of theirs.
 }
\end{remark}

\subsection{The case of $\JJ_-$}
We discuss briefly the case of $\JJ_-$. As mentioned previously, the situation corresponds to basically reversing the complex structure $I_-$. The Lie algebroid decomposition is
$$L_- = \ell_+ \dsum \bar \ell_- \Longrightarrow \wedge^*\bar L_-^* \cong \Dsum_{p,q} \wedge^p\bar\ell_+^* \tensor \wedge^q \ell_-^*$$
which gives the isomorphism as linear spaces
$$\Omega^*(\bar L_-) \cong \Dsum_{p,q} \Omega^{p,q}(\bar L_-), \text{ where } \Omega^{p,q}(\bar L_-) = \Omega^p(\bar\ell_+) \tensor \Omega^q(\ell_-)$$
Note that $\Omega^q(\ell_-) \cong \Omega^{q,0}_-(M)$. Instead of \eqref{eq:deltaoperators}, we have
\begin{equation}\label{eq:deltaprimeoperators}
 \bar\delta_+ := \bar\partial_+ + \bar\partial^+ : \Omega^{p,q}(\bar L_-) \to \Omega^{p+1,q}(\bar L_-) \text{ and } \delta_- := \partial_- + \partial^- : \Omega^{p,q}(\bar L_-) \to \Omega^{p,q+1}(\bar L_-)
\end{equation}
where 
$$\delta_\mp : \Omega^*(M) \xto{\nabla^\pm} \Omega^*(M) \tensor \Omega^1(M) \to \Omega^*(M) \tensor \Omega_\mp^{1,0}(M)$$
It then follows that $(\Omega^{p,q}(\bar L_-), \bar\delta'_+, \delta'_-)$ is a double complex and the associated spectral sequences computes $H^*(\bar L_-)$.

We obtain an $I_+$-holomorphic Lie subalgebroid $\mc B_+ \subset \mc E_+$ from $\bar L_-$, via the reduction with respect to $\bar\ell_+$. The underlying smooth Lie algebroid is given by $\ell_-$. The double complex $(\Omega^{p,q}(\bar L_-), \bar\delta'_+, \delta'_-)$ is a resolution of the holomorphic de Rham complex of $\mc B_+$ and thus we get
$$\HH^*(\mc B_+) \cong H^*(\bar L_-)$$
Similarly, a $\JJ_-$-holomorphic vector bundle $V$ corresponds to a $\mc B_+$-module, and we have
$$\HH^*(V; \mc B_+) \cong H^*(V; \JJ_-)$$

In the counterpart for $I_-$, we may either work with $-I_-$ or $-\JJ_-$. To stay in $\mc E_-$, we work with $-\JJ_-$ and the corresponding Lie algebroid $L_-$. Then, $L_-$ induces an $I_-$-holomorphic Lie subalgebroid $\mc B_- \subset \mc E_-$ via reduction with respect to $\bar \ell_-$. We have the similar double complex $(\Omega^{p,q}(L_-), \delta''_+, \bar\delta''_-)$ and the isomorphism of cohomology groups
$$\HH^*(\mc B_-) \cong H^*(L_-) \cong \bar{H^*(\bar L_-)}$$
For a $\JJ_-$-holomorphic vector bundle $V$, $\bar V$ is now naturally a $-\JJ_-$-holomorphic vector bundle and corresponds to a $\mc B_-$-module. We have
$$\HH^*(\bar V; \mc B_-) \cong H^*(\bar V; -\JJ_-) \cong \bar{H^*(V; \JJ_-)}$$
Similar to Proposition \ref{prop:isomorphiccategories}, the category of $\JJ_-$-holomorphic vector bundles is isomorphic to that of the $\mc B_+$-modules and $\CC$-antilinearly isomorphic to that of the $\mc B_-$-modules.

\subsection{Hodge decomposition}\label{subsect:genHodgedecomposition}
Recall that in the generalized Hodge decomposition associated to a generalized K\"ahler structure, $U_{-n,0} = C^\infty(U)$, where $U$ is the canonical line bundle of $\JJ_+$ (cf. \S \ref{subsect:spinorbundle}). We describe how the double complex $(\Omega^{p,q}(U; \bar L_+), D_+, D_-)$ for $U$ is related to the generalized Hodge decomposition.

Again, let $L:=L_+$ and $\bar L:= \bar L_+$. The pairing $2\<,\>$ defines the isomorphism $L \cong \bar L^*$:
 $$\mf X \mapsto 2\<\mf X, \bullet\> \text{ for } \mf X \in L$$
which extends to isomorphisms $\wedge^*L \cong \wedge^*\bar L^*$. Let $\mc D$ denote the differential on $C^\infty(\wedge^*L)$ induced by this identification, then for $\mf X \in C^\infty(L)$ and $s, t \in C^\infty(\bar L)$:
\begin{equation}\label{eq:gerstenhaberdifferential}
 (\mc D \mf X)(s, t) = 2(a(s)\<\mf X, t\> - a(t)\<\mf X, s\> - \<\mf X, s*t\>) = - 2\<\mf X * s, t\>
\end{equation}

\begin{lemma}\label{lemma:spinorCartanformula}
 For $\mf X \in C^\infty(L)$ and $\chi \in C^\infty(U)$, we have $\bar\partial_+^\gamma(\mf X\cdot \chi) + \mf X \cdot \bar\partial_+^\gamma \chi = \mc D\mf X \cdot \chi$.
\end{lemma}
{\it Proof:}
 By the decomposition of $d_\gamma$ in \eqref{eq:dgammadecomposition}, for $s, t \in C^\infty(\bar L)$, we compute
 \begin{equation*}
  \begin{split}
   & (t\wedge s)\cdot\left(\bar\partial_+^\gamma(\mf X\cdot \chi) + \mf X \cdot \bar\partial_+^\gamma \chi\right) = (t\wedge s)\cdot\left(d_\gamma(\mf X\cdot \chi) + \mf X \cdot d_\gamma \chi\right) \\
   = & (t\wedge s)\cdot \LLC_{\mf X}\chi = - t\cdot\left((\mf X * s) \cdot \chi\right) = - 2\<\mf X * s, t\>\chi = (\mc D\mf X)(s, t) \chi
  \end{split}
 \end{equation*}
 Here $\LLC_\mf X$ denotes the extended Lie derivative in the sense of Hu-Uribe \cite{HuUribe06}. Because $s\cdot \chi = t \cdot \chi = 0$, we have $(\mc D\mf X)(s, t) \chi = (t\wedge s)\cdot (\mc D\mf X) \cdot \chi$.
\qed

\begin{prop}\label{prop:twistedpartialdifferential}
 For $\alpha \in C^\infty(\wedge^k L)$ and $\chi \in C^\infty(U)$, we have
 $$\bar\partial_+^\gamma(\alpha \cdot \chi) = \mc D\alpha \cdot \chi + (-1)^k \alpha \cdot d_\gamma \chi$$
\end{prop}
{\it Proof:}
 Lemma \ref{lemma:spinorCartanformula} gives the case of $k = 1$. Suppose that the statement holds for $\alpha \in C^\infty(\wedge^k L)$. Let $\mf X \in C^\infty(L)$ and we have
 \begin{equation*}
  \begin{split}
   & d_\gamma(\mf X \cdot \alpha \cdot \chi) = \LLC_\mf X (\alpha\cdot \chi) - \mf X \cdot d_\gamma(\alpha \cdot \chi)
  \end{split}
 \end{equation*}
 Taking the component in $C^\infty(\wedge^{k+2} L) \cdot C^\infty(U)$, we obtain
 \begin{equation*}
  \begin{split}
   & \bar\partial_+^\gamma\left((\mf X \wedge \alpha) \cdot \chi\right) = \alpha\cdot \LLC_\mf X \chi - \mf X \cdot \bar\partial_+^\gamma(\alpha \cdot \chi) \\
   = & \alpha \cdot \mc D \mf X \cdot \chi - \mf X \cdot \left(\mc D\alpha \cdot \chi + (-1)^k \alpha \cdot d_\gamma \chi \right) \\
   = & \mc D(\mf X \wedge \alpha) \cdot \chi + (-1)^{k+1} (\mf X \wedge \alpha) \cdot d_\gamma \chi
  \end{split}
 \end{equation*}
 The statement then follows by mathematical induction.
\qed

\begin{theorem}\label{thm:spectralsequenceaspartofHodge}
 $(\Omega^{p,q}(U; \bar L), D_+, D_-) \cong (U_{r,s}, \bar\delta_+^\gamma, \bar\delta_-^\gamma)$, where $r = p+q - n$ and $s = p-q$.
\end{theorem}
{\it Proof:}
 The groups are identified by the pairing $2 \<,\>$ on $\TT M$ and the Clifford action:
 $$\Omega^{p,q}(U; \bar L) = C^\infty(\wedge^p\bar\ell_+^* \tensor \wedge^q \bar\ell_-^* \tensor U) \cong C^\infty(\wedge^p \ell_+ \tensor \wedge^q \ell_-) \tensor U_{-n, 0} \xto{\cong} U_{r,s}$$
 The identification of the differentials follows from Proposition \ref{prop:twistedpartialdifferential}.
\qed

The generalized Hodge decomposition for $H^*_\gamma(M)$ in \cite{Gualtieri0409} then implies the following.
\begin{corollary}\label{coro:canonicalcomputestwistedcoho}
 $H^*(U; \JJ_+) \cong H^*_\gamma(M)$.
 \qed
\end{corollary}

\section{Courant trivialization}\label{sect:Couranttrivialization}
Generalized K\"ahler geometry on $K$ can be described explicitly using the \emph{Courant trivialization} of $\TT G := TG \dsum T^*G$ given in  \cite{AlekseevBursztynMeinrenken} (Remark $3.4$, \emph{loc.} \emph{cit.}), where $G = K^\CC$ is the complexified Lie group. We recall the relevant constructions here.

\subsection{Notations}\label{subsect:Liegroupnotations}
Let $\mf g$ be the Lie algebra of $G$, which is identified with the \emph{right invariant} vector fields, then $\mf g = \mf k \tensor \CC$. For $a \in \mf g$, let $X_a^r$ denote the right invariant vector field with value $a$ at the identity $e \in G$, and respectively $\theta^r$ be the $\mf g$-valued right invariant Cartan $1$-form, i.e.
$$\theta^r(X^r_a) = a$$
The left invariant vector fields are denoted $X_a^l$ for $a \in \mf g$ and the left Cartan $1$-form is denoted $\theta^l$. The adjoint action relates the left and right invariant objects, for example, at $g \in G$
$$X_a^l(g) = X_{\Ad_g a}^r(g) \text{ with } \Ad_g a =  L_{g*} R_{g*}^{-1} a$$

Let $F:\mf g^{\tensor k} \to R$ be a $\CC$-linear map to a $\CC$-linear space $R$, then $F(\theta^r)$ is the \emph{$R$-valued} right invariant $k$-form on $G$, defined by skew-symmetrization of $F$:
$$(F(\theta^r))(X_{a_1^r}, \ldots, X_{a_p}^r) = \frac{1}{p!}\sum_{\sigma \in S_p} (-1)^\sigma F(a_{\sigma(1)}, \ldots, a_{\sigma(p)})$$
Similarly, the left-invariant form $F(\theta^l)$ is defined.
When $G$ is semi-simple, the Killing form $\kappa$ on $\mf g$ defines a bi-invariant Riemannian metric on $K$:
$$\sigma(X_a^r, X_b^r) = \kappa(a, b) = -\tr(\ad_a \circ \ad_b)$$
and the bi-invariant Cartan $3$-form $\gamma$ is given by $\gamma = \kappa(\theta^r, [\theta^r, \theta^r])$, or more explicitly
$$\gamma(X_a^r, X_b^r, X_c^r) = \kappa(a, [b, c]) =: \Lambda(a, b, c)$$
In general, $G = G' \times T'$, where $G'$ is semi-simple and $T'$ is an abelian Lie group. A bi-invariant metric can be defined on $K$ using the Killing form on $\mf g'$ and a non-degenerate bi-invariant form on $\mf t'$. This is the case for Hopf surfaces mentioned previously.

Let $\TT G = TG \dsum T^*G$, we have the Dorfman bracket defined by $\gamma$
$$(X + \xi) * (Y + \eta) = [X, Y] + \LLC_X \eta - \iota_Y d\xi + \iota_X\iota_Y \gamma$$
which defines the structure of Courant algebroid on $\TT G$. The restriction to $\TT K$ defines the generalized tangent bundle of $K$.

\subsection{Linear algebra}\label{subsect:spinorlinearalgebra}
Consider the isomorphism of linear spaces with quadratic forms:
\begin{equation}\label{eq:quadraticisomorphism}
 \kappa_* : \mf d = (\mf g \dsum \mf g, -\kappa \dsum \kappa) \to (\mf g \dsum \mf g^*, \<,\>) : (a, a') \mapsto (a'-a, \kappa(a) + \kappa(a'))
\end{equation}
which induces the isomorphism of Clifford algebras $\kappa_* : \Cl(\mf d) \xto\cong \Cl(\mf g \dsum \mf g^*)$. Write $\hat a := (a, a') \in \bar{\mf d}$, the natural spinor modules for these two Clifford algebras are
\begin{enumerate}
 \item $\wedge^* \mf g$ as $\Cl(\mf g \dsum \mf g^*)$-module: $(a, \xi) \circ v = a\wedge v + \iota_\xi v$ for $v \in \wedge^*\mf g$
 \item $\Cl(\mf g, \kappa)$ as $\Cl(\mf d)$-module: $\hat a \circ u = a' \cdot u - (-1)^{|u|} u \cdot a$ for $u \in \Cl(\mf g, \kappa)$
\end{enumerate}
The \emph{quantization map} $q : \wedge^* \mf g \to \Cl(\mf g, \kappa)$ is the unique spinor module (iso)morphism with respect to $\kappa_*$ such that $q(1) = 1$. By definition,
\begin{equation}\label{eq:spinormoduleisomorphismuptokappa}
 q(\kappa_*\hat a\circ v) = \hat a \circ q(v) \text{ for } v \in \wedge^*\mf g
\end{equation} 
which implies that $q|_{\mf g} = \id_\mf g$. Let $v \in \wedge^*\mf g$, then $q(v) \in \Cl(\mf g, \kappa)$. The left and right Clifford multiplications in $\Cl(\mf g, \kappa)$ interact with $q$ as following
$$a\cdot q(v) = q(a \wedge v + \iota_{\kappa(a)}v) \text{ and } (-1)^{|v|}q(v) \cdot a = q(a \wedge v - \iota_{\kappa(a)}v)$$

Switching $\mf g$ and $\mf g^*$, we have the $\Cl(\mf g \dsum \mf g^*)$-module $\wedge^*\mf g^*$:
$$(a, \xi) \circ \alpha = \iota_a \alpha + \xi \wedge \alpha \text{ for } \alpha \in \wedge^*\mf g^*$$
Fix a nonzero element $\mu \in \wedge^{2n} \mf g$. The \emph{star map} $\star: \wedge^*\mf g^* \to \wedge^{2n-*} \mf g$ is the unique spinor module (iso)morphism with $\star 1 = \mu$. It follows that for $a \in \mf g$ and $\xi \in \mf g^*$, we have $\star\xi = \iota_\xi \mu$ and
$$\star(\iota_a \alpha + \xi \wedge \alpha) = a \wedge \star \alpha + \iota_\xi (\star \alpha) \text{ for } \alpha \in \wedge^*\mf g^*$$
We thus have the isomorphism of spinor modules:
$$q \circ \star : \wedge^*\mf g^* \xto\cong \Cl(\mf g, \kappa)$$

Since the adjoint action of $G$ preserves $\kappa$, it extends to the adjoint action on $\Cl(\mf g, \kappa)$. We lift the adjoint action $\Ad : G \to SO(\mf g, \kappa)$ to a group homomorphism $\tau: G \to \Spin(\mf g, \kappa) \subset \Cl(\mf g, \kappa)$, and write $\tau_g = \tau(g)$, then
$$\Ad_g(u) = \tau_g \cdot u \cdot \tau_g^{-1} \text{ for } u \in \Cl(\mf g, \kappa)$$
This is possible because of $G$ is connected and $\pi_1(G)$ is torsionless.

\subsection{The trivializations}\label{subsect:thetrivializations}
The left and right actions of $G$ on itself define a left action of the \emph{double group} $D = G \times G$ on $G$:
$$\hat g \circ h = g' h g^{-1} \text{ for } \hat g = (g, g') \in D$$
The adjoint action of $G$ is induced by the diagonal embedding $G \into D$. Let $\mf d = \mf g \dsum \mf g$ be the Lie algebra, then the infinitesimal action of $\hat a = (a, a') \in \mf d$ is given by
$$X_{\hat a} = X_{a'}^r - X_a^l$$
It lifts to an extended action of $\mf d$ on $\TT G$
$$\mf d \to C^\infty(\TT G) : \hat a \mapsto \mf X_{\hat a} = X_{a'}^r - X_a^l + \kappa(\theta^r, a') + \kappa(\theta^l, a)$$
in the sense that the extended symmetry (cf. \cite{HuUribe06}) generated by $\mf X_{\hat a}$ is precisely $X_{\hat a}$. Let $\hat a = (a,a')$ and $\hat b = (b,b') \in \mf d$, direct computation shows
$$\<\mf X_{\hat a}, \mf X_{\hat b}\> = \kappa(a', b') - \kappa(a, b) \text{ and } \mf X_{\hat a} * \mf X_{\hat b} = \mf X_{[\hat a,\hat b]}$$
Thus we obtain the \emph{Courant trivialization} (cf. \cite{AlekseevBursztynMeinrenken})
$$G \times \mf d \xto{\cong} \TT G : \hat a \mapsto \mf X_{\hat a}$$
which coincides with $\kappa_*$ in \eqref{eq:quadraticisomorphism} when restricted to the identity $e \in G$.
The following will be useful.
\begin{lemma}[\cite{AlekseevBursztynMeinrenken}]\label{lemma:isotropicLiesubalgebra}
 Suppose $\mf l \subset \mf d$ is an isotropic Lie subalgebra, then the subbundle $L$ generated by $\mf X_{\hat a}$ with $\hat a \in \mf l$ is a Dirac structure, in particular is a Lie algebroid under the restriction of $*$.
 \qed
\end{lemma}

The isomorphism of Clifford algebras induced by $\kappa_*$ extends to an isomorphism of Clifford bundles
\begin{equation}\label{eq:cliffordisomorphismbykappa}
 G \times \Cl(\mf d) \to \Cl(TG \dsum TG, -\sigma \dsum \sigma) \xto{\kappa_*} \Cl(\TT G) : (g, \hat a) \mapsto (X_a^l, X_{a'}^r)(g) \mapsto \mf X_{\hat a}(g)
\end{equation}
Since $\mu \in \wedge^{2n} \mf g$ defines naturally a bi-invariant top multivector $\mu_G \in \Omega^*(G)$, we extend $\star$ pointwisely to $\star: \wedge^*T^*G \to \wedge^*TG$. Define $\Ad_{g*} = L_{g*} R_{g*}^{-1}$ on $\wedge^*T^*(G)$, then
$$\star \circ \Ad_{g} = \Ad_{g*} \circ \star : \wedge^* T^*G \to \wedge^* TG$$
The quantization map $q$ induces the isomorphism
\begin{equation}\label{eq:spinorisomorphismbyq}
 G \times \Cl(\mf g, \kappa) \to \Cl(TG, \sigma) \xto{q^{-1}} \wedge^*TG : (g, u) \mapsto X_{u\cdot \tau_g^{-1}}^r(g)  \mapsto q^{-1}(X_{u\cdot \tau_g^{-1}}^r(g) )
\end{equation}
where for $a_i \in \mf g$, $u = a_1 \cdot \ldots \cdot a_k \in \Cl(\mf g, \kappa)$ we write $X_u^r := X_{a_1}^r \cdot \ldots \cdot X_{a_k}^r$.
The isomorphism of spinor modules \eqref{eq:spinormoduleisomorphismuptokappa} implies that the Clifford action of terms in \eqref{eq:cliffordisomorphismbykappa} on the respective term in \eqref{eq:spinorisomorphismbyq} can be identified, which gives the isomorphism of the spinor modules.
\begin{lemma}[\cite{AlekseevBursztynMeinrenken}, Proposition $4.2$ a)]\label{lemma:spinortrivialization}
 The map $\mc Q : G \times \Cl(\mf g, \kappa) \to \wedge^*T^*G$ : 
 $$(g, u) \mapsto (q\circ \star)^{-1}\left(X_{u\cdot \tau_g^{-1}}^r(g)\right) \text{ for } u \in \Cl(\mf g, \kappa)$$
 is an isomorphism of spinor modules with respect to the isomorphism
 in \eqref{eq:cliffordisomorphismbykappa}.
 \qed
\end{lemma}

\subsection{Differentials} \label{subsect:spinordifferentials}
The adjoint action of $G$ on $\Cl(\mf g, \kappa)$ can be extended to a $D$-action:
$$\hat g \circ u = \tau_{g'} \cdot u \cdot \tau_g^{-1}$$
On the other hand, the adjoint action of $D$ on $\mf d$ induces a $D$-action on $\Cl(\mf d, -\kappa \dsum \kappa)$:
$$\Ad_{\hat g}(a, a') = (\Ad_g (a), \Ad_{g'} (a')) = \left(\tau_g \cdot a \cdot \tau_g^{-1}, \tau_{g'} \cdot a' \cdot \tau_{g'}^{-1}\right)$$
Then these two actions are compatible:
\begin{equation}\label{eq:compatibleactiononCliffordalgebras}
 \hat g \circ (\hat a \circ u) = \Ad_{\hat g}(\hat a)\circ (\hat g \circ u)
\end{equation}

The $D$-action on $G$ induces the action on $\Omega^*(G)$ by push-forward:
$$\hat g \circ \alpha = L_{g'*}R_{g*}^{-1} \alpha$$
For $u \in \Cl(\mf g, \kappa)$, let $\alpha_u(g) := \mc Q(g, u) \in \Omega^*(G)$. Let $g' \in G$, then we compute
$$R_{h*} (\alpha_u(g)) = R_{h*}(q\circ \star)^{-1}(X_{u \cdot \tau_g^{-1}}^r(g)) = (q\circ \star)^{-1}(X_{u\cdot \tau_g^{-1}}^r(gh) ) = \alpha_{u\cdot \tau_{h}}(gh)$$
$$L_{h*} (\alpha_u(g)) = \Ad_{h*} (q\circ \star)^{-1}(X_{u\cdot \tau_g^{-1}}^r(gh) ) = (q\circ \star)^{-1}(X_{\tau_{h} \cdot u\cdot \tau_g^{-1} \cdot \tau_{h}^{-1}}^r(hg)) = \alpha_{\tau_{h}\cdot u}(hg)$$
Thus, $\mc Q$ is $D$-equivariant (\cite{AlekseevBursztynMeinrenken} Proposition $4.2$ c)), i.e.
$$\hat g\circ \alpha_u = \alpha_{\hat g \circ u} \text{ for } u \in \Cl(\mf g, \kappa)$$

The homomorphism $\tau$ induces an homomorphism $\tau' : \mf g \to \Cl(\mf g, \kappa)$:
$$\tau'_a := \left.\frac{d}{dt}\right|_{t = 0} \tau_{\exp(ta)} \text{ for } a \in \mf g$$
which implies that
$$[a, b] = \ad_a b = \tau'_a \cdot b - b \cdot \tau'_a \text{ and }\tau'_{[a,b]} = [\tau'_a, \tau'_b]$$
Let $\{e_i\}_{i=1}^{2n}$ and $\{e^i\}_{i=1}^{2n}$ be dual basis of $\mf g$, e.g. a $\kappa$-orthonormal $\RR$-basis $\{e_i = e^i\}_{i=1}^{2n}$ of $\mf k$, 
and let 
$$\Theta := q(\Lambda) = \frac{1}{6}\sum_{i,j,k = 1}^{2n}\Lambda(e_i, e_j, e_k) e^i\cdot e^j \cdot e^k$$ 
then direct computation shows (\cite{AlekseevBursztynMeinrenken}, Meinrenken \cite{Meinrenken13})
\begin{equation}\label{eq:Cliffordadjoint}
 \tau'_a = \frac{1}{4} \sum_{i= 1}^{2n} [a, e_i] \cdot e^i = -\frac{1}{4} [\Theta, a] \in \Cl(\mf g, \kappa)
\end{equation}
The infinitesimal action of $\hat a = (a, a') \in \mf d$ is given by the Clifford action of $\tau'_{\hat a} := (\tau'_{a}, \tau'_{a'}) \in \Cl(\bar{\mf d})$
$$\tau'_{\hat a} \circ u = \tau'_{a'} \cdot u - u \cdot \tau'_a \text{ for } u \in \Cl(\mf g, \kappa)$$
since $\tau'_a$ is even.
The infinitesimal action of $\hat a \in \mf d$ on $\Omega^*(G)$ is given by \footnote{We use push-forward to define the Lie derivative here.}
$$\LLC_{X_{\hat a}} \alpha = \left.\frac{d}{dt}\right|_{t = 0} \exp(t\hat a) \circ \alpha$$
It then follows from $D$-equivariance of $\mc Q$ that
\begin{equation}\label{eq:infinitesimalinvariance}
 \alpha_{\tau'_{\hat a} \circ u} = \LLC_{X_{\hat a}}\alpha_u
\end{equation}

From \cite{AlekseevBursztynMeinrenken}, $\displaystyle{d^\Cl := \frac{1}{4}[\Theta, \bullet]}$ is the \emph{Clifford differential} on $\Cl(\mf g, \kappa)$ (cf. \cite{Meinrenken13} Chap. $6$), i.e. $d^\Cl \circ d^\Cl = 0$.

\begin{lemma}[\cite{AlekseevBursztynMeinrenken} Proposition $4.2$ b)]\label{lemma:differentialtrivialization}
 $\alpha_{d^\Cl u} = d_\gamma \alpha_u$.
\end{lemma}
{\it Proof:}
 Using Clifford action of $C^\infty(\TT G)$ on $\Omega^*(G)$, we have (cf. \cite{HuUribe06})
 $$\LLC_{X_{\hat a}} \alpha = -\mf X_{\hat a} \circ d_\gamma \alpha - d_\gamma(\mf X_{\hat a} \circ \alpha)$$
 where $(X + \xi) \circ \alpha = \iota_X \alpha + \xi \wedge \alpha$ is the Clifford action of $C^\infty(\TT G)$ on $\Omega^*(G)$.
 For the infinitesimal action on $\Cl(\mf g, \kappa)$, direct computation shows that
 $$\tau'_{\hat a} \circ u = \tau'_{a'}\cdot u - u \cdot \tau'_a = -\hat a \circ d^\Cl u - d^\Cl(\hat a \circ u)$$
 The statement follows from \eqref{eq:infinitesimalinvariance} and that $d^\Cl$ ($d_\gamma$) and $\mc Q$ have opposite parities (cf. \cite{AlekseevBursztynMeinrenken}).
\qed

\section{Lie algebraic generalized K\"ahler structures} \label{sect:canonicalgenkalher}

{\bf Notation:} The notation $\mf l \bdsum \mf l' \subseteq \mf g \dsum \mf g'$ means that $\mf l \subseteq \mf g$ and $\mf l' \subseteq \mf g'$, while $\mf p \dsum \mf p' \subseteq \mf V$ means that $\mf p$ and $\mf p'$ are both subspaces of $\mf V$.

Let $2n = \dim_\RR \mf k$, then $2n = \dim_\CC \mf g$. 
Let $\mf t$ be a Cartan subalgebra of $\mf k$ then $\mf h = \mf t \tensor_\RR \CC$ is a Cartan subalgebra of $\mf g$. Let $R$ be the set of roots with respect to $\mf h$. Since $\mf k$ is a compact Lie algebra, any root $\alpha \in R$ is purely imaginary on $\mf t$. In particular, when $\mf k$ is semisimple, the Satake diagram of $\mf k$ is simply the Dynkin diagram of $\mf g$ with all vertices painted black.

Let $\mf b$ be a Borel subalgebra containing $\mf h$ and $R_+ \subset R$ (respectively $R_-$) be the corresponding set of positive (respectively negative) roots. We then have the triangular decomposition
$$\mf g = \mf n_+ \dsum \mf h \dsum \mf n_- \text{ where } \mf n_\pm = \dsum_{\alpha \in R_\pm} \mf g_\alpha$$
When $\mf k$ is semisimple, this is an orthogonal decomposition with respect to $\kappa$. 
In general, let $\mf z$ be the center of $\mf g$ then we can decompose $\mf h$ furthermore as 
$$\mf h = \mf h_0 \dsum \mf z \text{ where } \mf h_0 = \mf h \inter [\mf g, \mf g]$$
It then follows that $\mf g = [\mf g, \mf g] \dsum \mf z$, and we extend $\kappa$ from the semisimple part $[\mf g, \mf g]$ to $\mf g$ such that this decomposition is orthogonal.

There are root vectors $a_\alpha \in \mf g_\alpha$ for $\alpha \in R$ (cf. Humphreys \cite{Humphreys}) such that
$$[a_\alpha, a_\beta] = a_{\alpha + \beta} \text{ for } \alpha \pm \beta \neq 0, \text{ and } \kappa(a_\alpha, a_{-\alpha}) = 1$$
Since $[h, a_\alpha] = \alpha(h) a_\alpha$ for any $h \in \mf t$, we see that $\mf g_\alpha \inter \mf k = \{0\}$. Let $h_\alpha = [a_\alpha, a_{-\alpha}]$ for $\alpha \in R_+$, then $\kappa(h_\alpha, h) = \alpha(h)$ for all $h \in \mf h$. It follows that
$$h_\alpha \in i\mf t \subset \mf h \text{ for all } \alpha \in R$$

\subsection{Lagrangian subalgebras}\label{subsect:Lagrangiansubalgebra}
The classification of complex Lagrangian subalgebra of in $\mf d$ (with respect to $-\kappa \dsum \kappa$) is obtained by Karolinsky \cite{Karolinsky00}. We recall the description of 
a representative for each orbit of adjoint action of $D$ on the space of Lagrangian subalgebras of $\mf d$, as given in Evens-Lu \cite{EvensLu06}. Let $\mf b \supset \mf h$ be a Borel subalgebra of $\mf g$ and $R_+$ the corresponding set of positive roots. Let $\Gamma \subset R_+$ be the set of simple roots. For any $P \subseteq \Gamma$, we set
$$[P] := R \inter \Span_\CC\{\alpha \in P\}$$
Consider the corresponding subalgebra $\mf g_P \subset \mf g$ defined by the sub-diagram of the Dynkin diagram of $\mf g$ with vertices in $P$, then $\mf m_P := \mf g_P + \mf h$ decomposes as $\mf m_P = \mf g_P \dsum \mf z_P$, where $\mf z_P$ is the center of $\mf m_P$. We have also the nilpotent subalgebras
$$\mf n_P := \dsum_{\alpha \in R_+ \setminus [P]} \mf g_\alpha \text{ and }\mf n_P^- := \dsum_{\alpha \in R_+ \setminus [P]} \mf g_{-\alpha}$$
A \emph{generalized Belavin-Drinfeld triple} $(P, P', \pi)$ consists of two subsets $P,P' \subseteq \Gamma$ and an \emph{isometry} $\pi : P \to P'$, i.e. an isomorphism such that
$$\kappa(h_\alpha, h_\beta) = \kappa(h_{\pi\alpha}, h_{\pi\beta})$$
It induces a unique Lie algebra isomorphism $\psi_\pi : \mf g_P \to \mf g_{P'}$ such that
$$\psi_\pi(a_\alpha) = a_{\pi\alpha} \text{ for all } \alpha \in P$$
\begin{lemma}[\cite{EvensLu06}, Theorem $2.16$]\label{lemma:EvensLurepresentative}
 Each orbit of the adjoint action of $D$ on the space of Lagrangian subalgebras of $\mf d$ contains exactly one element of the form
 $$\mf l(\pi, F) := F \dsum (\mf n_P \dsum \mf n_{P'}^-) \dsum L^\pi$$
 where $(P,P', \pi)$ is a generalized Belavin-Drinfeld triple, $F \subset \mf z_P \dsum \mf z_{P'}$ is a Lagrangian subspace, and $L^\pi$ is the graph of $\psi_\pi$.
 \qed
\end{lemma}

\subsection{Invariant complex structures} \label{subsect:invariantcomplexstructures}
Since $\mf g = \mf k \tensor_\RR \CC$, a complex structure $J$ on $\mf k$ corresponds to an $n$-dimensional complex subspace $\mf l_J \subset \mf g$ such that $\mf l_J \inter \mf k = \{0\}$. We say that $J$ is \emph{integrable} if $\mf l_J$ is a complex Lie subalgebra, called a \emph{Samelson subalgebra} (cf. Samelson \cite{Samelson53}). In this case, the corresponding invariant almost complex structures on $K$ with value $J$ at $e \in K$ are integrable.

Pittie \cite{Pittie88} showed that there is a unique Borel subalgebra $\mf b \subseteq \mf g$ such that $\mf l_J \subseteq \mf b$. Furthermore, an integrable complex structure $J$ on $\mf k$ is in $1$ -- $1$ correspondence to the pairs $(\mf b, J')$ where
\begin{enumerate}
 \item $\mf b \subset \mf g$ is a Borel subalgebra, and $\mf h' \in \mf b$ is the corresponding complex Cartan subalgebra
 \item $J'$ is a complex structure on the corresponding real Cartan subalgebra $\mf t' = \mf h' \inter \mf k$
\end{enumerate}
We say that $\mf b$ or $\mf t'$ are the Borel subalgebra or real Cartan subalgebra \emph{corresponding to} the complex structure $J$.
\begin{prop}[\cite{Pittie88} Proposition $2.6$, Corollary $2.7.1$]\label{prop:spaceofcomplexstructures}
 The space of (left) invariant $\sigma$-orthogonal integrable complex structures on $K$ is isomorphic to $O(\mf h) / U(\mf h) \times K / T$, and the moduli space up to $\Aut(K)$ is isomorphic to $F \backslash O(\mf h) / U(\mf h)$, where $F$ is the outer automorphism group of $K$.
 \qed
\end{prop}

Let $J$ be an integrable complex structure on $\mf k$ and let $\mf t_{1,0} \subset \mf h$ be the $i$-eigensubspace, then
$$\mf l_J = \mf n_+ \dsum \mf t_{1,0} \text{ and } \bar{\mf l}_J = \mf n_- \dsum \mf t_{0,1}$$
We note that the complex structure $J$ on $\mf k$ can also be induced from the exact sequence
\begin{equation}\label{eq:Liealgebracomplexstructure}
 0 \to \bmf l_J \to \mf g \to \mf k \to 0
\end{equation}
From now on, we assume that choices are made such that $J$ is orthogonal with respect to $\kappa$ on $\mf k$. It implies that $\mf l_J$ as well as $\bmf l_J$ are isotropic subspaces of $\mf g$, with respect to $\kappa$.

Let $J_\pm$ be two integrable complex structures on $\mf k$. We denote by $I_+$ (respectively $I_-$) the right invariant (respectively left invariant) complex structure on $K$ defined by $J_+$ (respectively $J_-$), e.g.
$$I_+(X_a^r) = X_{J_+a}^r \text{ and } I_-(X_a^l) = X_{J_-a}^l$$
It follows that the holomorphic tangent bundle of $I_+$ (respectively $I_-$) is spanned by $X_a^r$ (respectively $X_a^l$) for $a \in \mf l_+ := \mf l_{J_+}$ (respectively $a \in \mf l_- := \mf l_{J_-}$).
Unless $K$ is abelian, i.e. a torus, the right invariant complex structure $I_+$ is not left invariant. On the other hand, both $I_\pm$ are invariant under the action of the maximal torus generated by $\mf t$, since
$$[\mf t, \mf l_J] \subseteq \mf l_J$$

The following result will be useful later.
\begin{lemma}\label{lemma:alwaysGauduchon}
 $\sigma$ is Gauduchon with respect to both $I_+$ and $I_-$.
\end{lemma}
{\it Proof:}
 We verify the $+$-case. Let $\omega_+ = \sigma I_+$ be the K\"ahler form and $d\vol_\sigma$ denote the invariant volume form with respect to $\sigma$. Since both $I_+$ and $\omega_+$ are right-invariant, $dd_+^c\left(\omega_+^{n-1}\right)$ is a right-invariant top form, i.e. $dd_+^c\left(\omega_+^{n-1}\right) = A d\vol_\sigma$ for some $A \in \RR$. Since $K$ is compact
 $$0 = \int_K dd_+^c\left(\omega_+^{n-1}\right) = \int_K A d\vol_\sigma$$
 Thus $A = 0$ and $dd_+^c\left(\omega_+^{n-1}\right) = 0$, i.e. $\sigma$ is Gauduchon with respect to $I_+$.
\qed

It follows that for any $I_\pm$-holomorphic vector bundle $V$ on $K$, we have the notion of \emph{$\pm$-degrees} :
$$\deg_\pm (V) := \frac{i}{2} \frac{\displaystyle{\int_M F_\pm \wedge \omega_\pm^{n-1}}}{\displaystyle{\int_M \omega_\pm^{n-1}}}$$
where $F_\pm$ are the curvatures of the Chern connections with respect to any Hermitian metric on $V$. For $\alpha \in (0,1)$, we have the $\alpha$-degree (\cite{HuMoraruSeyyedali}): 
$$\deg_\alpha (V) = \alpha \deg_+ (V) + (1-\alpha) \deg_- V$$
The \emph{$\pm$-slope} $\mu_\pm$ and \emph{$\alpha$-slope} $\mu_\alpha$ are defined as quotients of the corresponding degrees by the rank.

\subsection{Bismut connections} \label{subsect:BismutconnectionsonLiegroup}
The Bismut connections $\nabla^\pm$ on $K$ for the Hermitian structures given by $I_\pm$ has respectively torsions $\pm \gamma$. Recall that $I_\pm$ is parallel in $\nabla^\pm$ respectively. They can also be defined from the Dorfman bracket as following. Let $X, Y \in TK$, then
$$\nabla^\pm_X Y = a\left\{[X \mp \sigma(X) * (Y \pm \sigma(Y))]^\pm\right\} = \nabla_X Y \pm \frac{1}{2} \sigma^{-1}\iota_X \iota_Y \gamma$$
where $a$ is the projection to $TK$ and $[\bullet]^\pm$ denote the projections to $C_\pm$.
\begin{lemma}\label{lemma:Bismutflat}
 $K$ is flat with respect to both $\nabla^\pm$, i.e. $TK$ can be trivialized by $\nabla^\pm$-horizontal sections.
\end{lemma}
{\it Proof:}
 We prove the $+$-case. For $a, b' \in \mf k$, let $d_1 = (a, 0)$ and $d_2 = (0, b') \in \mf d$, then $[d_1, d_2] = 0$. It follows that 
 $$\mf X_{d_1} * \mf X_{d_2} = \mf X_{[d_1, d_2]} = 0$$
 Since $\mf X_{d_1} = - X_a^l + \kappa(\theta^l, a) \in C_-$ and $\mf X_{d_2} = X_{b'}^r + \kappa(\theta^r, b') \in C_+$, it implies that $\nabla_{X_a^l}^+ X_{b'}^r = 0$. The result follows from the fact that the left invariant vector fields span $TK$.
\qed

\begin{corollary}\label{coro:formBismutflat}
 Suppose that $\alpha \in \Omega^p(K)$ is a right (respectively, left) invariant form, then $\nabla^+ \alpha = 0$ (respectively, $\nabla^- \alpha = 0$).
\end{corollary}
{\it Proof:}
 We prove the case of right invariant forms. Let $X_i$ be right invariant vector fields, for $i = 1, \ldots, p$, then for any vector field $X$ on $K$, we have
 $$X (\alpha(X_1, \ldots, X_p)) = (\nabla^+_{X}\alpha)(X_1, \ldots, X_p) + \sum_{i = 1}^p(-1)^p \alpha(\nabla^+_X X_i, X_1, \ldots, \hat X_i, \ldots, X_p)$$
 The lefthand side vanishes because right invariant functions are constants. The terms in the summation vanish because all the covariant derivatives vanish by Lemma \ref{lemma:Bismutflat}. Thus, $\nabla^+ \alpha = 0$.
\qed

\subsection{Generalized K\"ahler structures} \label{subsect:LiegroupgenKahler}
By Lemma \ref{lemma:EvensLurepresentative}, a complex Lagrangian subalgebra $\mf L \subset \mf d$ is of the form
\begin{equation}\label{eq:complexLagrangiangenform}
 \mf L = \Ad_{(g,g')}\mf l(\pi, F) = \Ad_{(g, g')}F \dsum (\Ad_g \mf n_P \dsum \Ad_{g'} \mf n_{P'}^-) \dsum \Ad_{(g,g')}L^\pi
\end{equation}
Lemma \ref{lemma:isotropicLiesubalgebra} implies that it defines a complex Dirac structure $L \subset \TT_\CC K$ via Courant trivialization. If furthermore $\mf L \inter \mf k \dsum \mf k = \{0\}$, the corresponding $L$ defines a generalized complex structure $\JJ$ on $K$. 
\begin{defn}\label{defn:Liealgebraicgencomplex}
 A generalized complex structure on $K$ defined via the Courant trivialization by a complex Lagrangian subalgebra $\mf L \subset \mf d$ satisfying $\mf L \inter (\mf k\dsum \mf k) = \{0\}$ is called a \emph{Lie algebraic generalized complex structure}. A generalized K\"ahler structure $(K, \gamma; \JJ_\pm, \sigma)$ is a \emph{Lie algebraic generalized K\"ahler structure} if both $\JJ_\pm$ are Lie algebraic generalized complex structures.
\end{defn}
\begin{lemma}\label{lemma:LiealgebraicgenKahler}
 Let $(K, \gamma; \JJ_\pm, \sigma)$ be a Lie algebraic generalized K\"ahler structure and $\mf L_\pm \subset \mf d$ be the corresponding complex Lagrangian subalgebra defining $\JJ_\pm$. Then there exist two (not necessarily distinct) Samelson subalgebras $\mf l_\pm \subset \mf g$ such that $\mf L_+ = \mf l_- \bdsum \mf l_+$ and $\mf L_- = \bmf l_- \bdsum \mf l_+$.
\end{lemma}
{\it Proof:}
 Let $L_\pm$ be the $i$-eigensubbundles of $\JJ_\pm$ respectively, then we have $L_+ = \ell_+ \dsum \ell_-$ and $L_- = \ell_+ \dsum \bar \ell_-$ (cf. \S\ref{subsect:genKahlergeomdescription}), where $\ell_\pm \subset C_\pm$. Via the Courant trivialization, by \eqref{eq:quadraticisomorphism}, $\ell_+$ (respectively $\ell_-$) is defined by a subspace $\mf l'_+$ of $0 \dsum \mf g$ (respectively a subspace $\mf l'_-$ of $\mf g \dsum 0$). Thus, we have, for example
 $$\mf L_+ = \mf l'_- \bdsum \mf l'_+ \subset \mf g \dsum \mf g$$
 In the general form \eqref{eq:complexLagrangiangenform}, since the adjoint action of $D$ does not affect splitness, this implies that
 \begin{enumerate}
  \item $L^\pi$ is trivial, i.e. the generalized Belavin-Drinfeld triple is trivial
  \item $\mf n_P$ (respectively $\mf n_{P'}^-$) is the sum of positive (respectively negative) root spaces
  \item $F = F_- \dsum F_+$ where $F_\pm \subset \mf h$ are complex Lagrangian subspaces of $\mf h$
 \end{enumerate}
 This proves the lemma.
\qed

Thus, the space of Lie algebraic generalized K\"ahler structures can be identified with the space of pairs of Samelson subalgebras, which admits the natural adjoint action by $K \times K$. Let $\mf b_\pm$ denote the unique Borel subalgebras that contain $\mf l_\pm$ respectively, then for any pair $(\mf l_+, \mf l_-)$ of Samelson subalgebras, there exists $(g,g') \in D$ such that $\Ad_{g'} \mf b_+ = \Ad_{g} \mf b_-$.
\begin{defn}\label{defn:inducedgenKahler}
 An \emph{induced} generalized K\"ahler structure on $K$ is a Lie algebraic generalized K\"ahler structure with $\mf b_+ = \mf b_-$. It is a \emph{canonical} generalized K\"ahler structure if $\mf l_+ = \mf l_-$. Two Lie algebraic generalized K\"ahler structures are \emph{algebraically equivalent} if the corresponding pairs of Samelson subalgebras are related by the adjoint action of $D$. They are \emph{geometrically equivalent} if the corresponding pairs of Samelson subalgebras are related by the adjoint action of $K \times K$.
\end{defn}

\begin{lemma}\label{lemma:algebraicequivalenttoinduced}
 Geometrically equivalent Lie algebraic generalized K\"ahler structures are isomorphic.
 Any Lie algebraic generalized K\"ahler structure on $K$ is algebraically equivalent to an induced generalized K\"ahler structure.
\end{lemma}
{\it Proof:}
 Let $(g',g) \in K \times K$, define $(\mf l_+, \mf l_-) = (\Ad_{g'} \mf l_+', \mf l_-)$. For the corresponding invariant complex structures, we have
 $$I_+ = \Ad_{g'*} I_+' = L_{g'*} I_+'$$
 Since $I_-$ is left invariant, we see that $(I_+, I_-) = L_{g'*}(I_+', I_-)$, i.e. $L_{g'}$ induces an isomorphism of the bi-Hermitian structures. Similarly, we see that when $(\mf l_+, \mf l_-) = (\mf l_+', \Ad_g\mf l_-)$, the isomorphism is induced by $R_{g}$. The second statement is straightforward.
\qed

The Lie algebraic generalized complex structures can be explicitly written down.
In terms of $I_\pm$, for $a, a' \in \mf k$, we have:
$$\JJ_+ \left(X_{a'}^r - X_{a}^l + \sigma(X_{a'}^r + X_{a}^l)\right) = I_+X_{a'}^r - I_-X_{a}^l + \sigma(I_+X_{a'}^r + I_-X_{a}^l) $$
$$\JJ_- \left(X_{a'}^r - X_{a}^l + \sigma(X_{a'}^r + X_{a}^l)\right) = I_+X_{a'}^r + I_-X_{a}^l + \sigma(I_+X_{a'}^r - I_-X_{a}^l)$$
Let
$$\ell_+ = \Span_\CC\{\mf X_{0,a} : a \in \mf l_+\} \text{ and } \ell_- = \Span_\CC\{\mf X_{a, 0} : a \in \mf l_-\}$$
then the $i$-eigensubbundles are
$$L_+ = \ell_+ \dsum \ell_- = \Span_\CC\{\mf X_{a, a'} : (a, a') \in \mf l_- \bdsum \mf l_+\} \text{ and } L_- = \ell_+ \dsum \bar\ell_- = \Span_\CC\{\mf X_{a, a'} : (a, a') \in \bmf l_- \bdsum \mf l_+\}$$

\begin{lemma}\label{lemma:kahlermaximaltorus}
 Let $(K, \gamma; \JJ_\pm, \sigma)$ be a canonical generalized K\"ahler structure. Let $T \subset K$ be the maximal torus generated by the real Cartan subalgebra $\mf t$ (corresponding to $J = J_+ = J_-$). Then, the restriction of $(K, \gamma; \JJ_\pm, \sigma)$ on $T$ is an invariant K\"ahler structure.
\end{lemma}
{\it Proof:}
 For $a \in \mf t$ and $h \in T$, we have $X_a^r(h) = X_a^l(h) =: X_a(h)$. Since $J$ preserves $\mf t$, we have $I_\pm X_a = X_{Ja}$. It follows that when restricted to $T$, we have
 $$\JJ_+(X_{a'-a} + \sigma(X_{a'+a})) = X_{J(a'-a)} + \sigma(X_{J(a'+a)}) \text{ and } \JJ_-(X_{a'-a} + \sigma(X_{a'+a})) = X_{J(a'+a)} + \sigma(X_{J(a'-a)})$$
 It's now clear that $\JJ_+$ is defined by the invariant complex structure $I := I_\pm |_T$ and $\JJ_-$ is defined by the invariant symplectic structure $\omega := \omega_\pm |_T$. They define an invariant K\"ahler structure.
\qed

As invariant complex structures are determined by the restriction to a maximal torus up to group isomorphisms (cf. \cite{Pittie88}), the same holds for canonical generalized K\"ahler structures.
\begin{corollary}\label{coro:torusdeterminegroup}
 Up to group automorphisms, the canonical generalized K\"ahler structures on $K$ are determined by the invariant K\"ahler structures on a maximal torus $T$.
 \qed
\end{corollary}

We may apply the construction directly to the even dimensional torus $T$ and obtain generalized K\"ahler structures on it. Since $T$ is abelian, the resulting generalized K\"ahler structures are invariant under $T$ action, which will be called the \emph{invariant generalized K\"ahler structures} on $T$.
Corollary \ref{coro:torusdeterminegroup} generalizes to the statement about induced generalized K\"ahler structures.
\begin{prop}\label{prop:torusdeterminegroup}
 Up to group automorphisms, induced generalized K\"ahler structures on $K$ are determined by invariant generalized K\"ahler structures on a maximal torus $T$.
 \qed
\end{prop}

Since the left and right invariant complex structures are invariant under the actions of their respective maximal tori generated by $\mf t_\pm \subseteq \mf b_\pm$, the resulting generalized K\"ahler structure is invariant under the action of a subtorus, i.e. the intersection. In this sense, induced structures carry the maximal amount of symmetries.
\begin{prop}\label{prop:inducedtorusinvariant}
 An induced generalized K\"ahler structure is invariant under the action of a maximal torus, which induces naturally on the quotient the homogeneous K\"ahler structure.
 \qed
\end{prop}

\subsection{Canonical bundles}\label{subsect:genKahlercanonicalbundle}
We have first the vanishing result for semi-simple Lie groups.
\begin{prop}\label{prop:canonicalcohotrivial}
 $H^*(U_\pm; \JJ_\pm) = 0$ when $K$ is a semi-simple Lie group.
\end{prop}
{\it Proof:}
 It's well-known that the twisted de Rham cohomology $H^*_\gamma(K)$ is identically trivial for semi-simple $K$, e.g. \cite{Ferreira13}. The result follows from Corollary \ref{coro:canonicalcomputestwistedcoho}.
\qed
 
In general, the $\JJ_\pm$-canonical bundles on $K$, using $\mc Q$ in Lemma \ref{lemma:spinortrivialization} (cf. \cite{AlekseevBursztynMeinrenken}), are determined by the corresponding spinor lines in $\Cl(\mf g, \kappa)$.
Let $\bmf p \subseteq \mf t_{-,0,1}$ be a subspace such that
$$\bmf m := \bmf l_+ + \bmf l_- = \bmf l_+ \dsum \bmf p \subseteq \mf g$$
Let $\{\bar b_j\}_{j = 1}^{n+s}$ be a basis of $\bmf m$ such that 
\begin{equation}\label{eq:basiscondition}
 \{\bar b_j\}_{j=1}^n \text{ is a basis of } \bmf l_+, \{\bar b_j\}_{j=s+1}^{n+s} \text{ is a basis of } \bmf l_-, \text{ and } \{\bar b_j\}_{j = n+1}^{n+s} \text{ is a basis of } \bmf p
\end{equation}
and define
\begin{equation}\label{eq:pluspinordefn}
 u_+ = \bar b_1 \cdot \ldots \cdot \bar b_{n+s} \in \Cl(\mf g, \kappa)
\end{equation}
Similarly we can define $u_-$ from $\bmf l_+ + \mf l_- \subseteq \mf g$.
\begin{lemma}\label{lemma:inducedstructureplusspinor}
  $u_+$ as given above is a pure spinor defining $\bmf l_- \bdsum \bmf l_+ \subset (\mf d, -\kappa \dsum \kappa)$, and the $\JJ_+$-canonical bundle $U_+$ is generated by $\chi_+ := \mc Q(u_+) \in \Omega^*(K)$. The type of $\JJ_+$ (respectively of $\JJ_-$) has the same parity as $\dim_\CC \bmf l_+ \inter \bmf l_-$ (respectively as $\dim_\CC \bmf l_+ \inter \mf l_-$).
\end{lemma}
{\it Proof:}
 Left as an exercise for the reader.
\qed

For an induced generalized K\"ahler structures, where a Borel subalgebra $\mf b$ contains both $\mf l_\pm$, the spinor $u_+$ admits a more explicit description. Let 
$$\displaystyle{\mf b = \dsum_{\alpha \in R_+} \mf g_\alpha \dsum \mf h}$$ 
be the decomposition of the Borel subalgebra $\mf b$ into the positive root spaces $\mf g_\alpha$ and the Cartan subalgebra $\mf h$. Correspondingly, 
$$\displaystyle{\bmf l_\pm = \dsum_{\alpha \in R_+} \mf g_{-\alpha} \dsum \mf t_{\pm, 0,1}} \text{ and } \bmf p \subseteq \mf h$$
Let $\{\bar h_i: \bar h_i \in \mf t_{+, 0,1}\}$ be a basis of $\mf t_{+, 0,1}$ and $\{\bar p_1, \ldots, \bar p_s\}$ be a basis for $\bmf p$. Choose $0 \neq a_{-\alpha} \in \mf g_{-\alpha}$, then
$$\{\bar b_j\} := \{\bar h_i, a_{-\alpha}, \bar p_k\} \text{ is a basis of } \bmf m$$
and satisfies the conditions in \eqref{eq:basiscondition} and defines correspondingly the spinor $u_+$ is this case.

By Lemma \ref{lemma:algebraicequivalenttoinduced} we have
\begin{corollary}\label{coro:generalstructurepluspinor}
 The $\JJ_+$-canonical bundle of a Lie algebraic generalized K\"ahler structure is of the form $\<\mc Q(\hat g \circ u_+)\>$, where $\hat g \in D$ and $u_+$ is a spinor for an induced structure, as described above.
\qed
\end{corollary}

The $\JJ_+$-holomorphic structure on $U_+$ is given by the twisted differential $d_\gamma$ on $C^\infty(U_+) \subset \Omega^*(M)$ (cf. \S \ref{subsect:spinorbundle}, \cite{Gualtieri04}), which corresponds to the Clifford differential $d^\Cl$ under the Courant trivialization (cf. Lemma \ref{lemma:differentialtrivialization}, \cite{AlekseevBursztynMeinrenken}). 
Let $\{b_j\}_{j=1}^n$ be the basis for $\mf l_ +$ dual to $\{\bar b_j\}_{i=1}^n$, i.e.
$$\{b_j\}_{i = 1}^n = \{h_j : h_j \in \mf t_{+,1,0}\} \union \{a_{\alpha} \in \mf g_{\alpha} : \alpha \in R_+\}$$
such that $\kappa(b_i, \bar b_j) = \delta_{ij}$. Let $h_\alpha = [a_\alpha, a_{-\alpha}] \in \mf h$, then we have $\kappa(h_\alpha, h) = \alpha(h)$ for all $h \in \mf h$.
$$u_{\mf l_+} = \bar b_1 \cdot \ldots \cdot \bar b_n \text{ and } u_{\mf p} = \bar p_1 \cdot \ldots \cdot \bar p_s$$
then $u_+ = u_{\mf l_+} \cdot u_{\mf p}$. We have the \emph{Weyl vector} \footnote{In Lie theory, $\rho$ is often used to denote the \emph{Weyl vector} in $\mf h^*$, i.e. half sum of positive roots. We use the same notation to denote the corresponding element in $\mf h$ via $\kappa$. In particular, $\rho \in i \mf t$.}
$$\displaystyle{\rho = \frac{1}{2} \sum_{\alpha \in R_+} h_\alpha \in \mf h}$$
We note that the norm of the Weyl vector with respect to $\kappa$ does not depend on the choice of the Borel subalgebra $\mf b$, and will be denoted $|\rho|$.
\begin{lemma}\label{lemma:ClifforddifferentialSamelson}
 $\displaystyle{d^\Cl u_{\mf l_+} = \frac{1}{2}(-\rho, \rho) \circ u_{\mf l_+}}$.
\end{lemma}

{\it Proof:}
 Let $\mf l:=\mf l_+$.
 In the basis $\{b_j, \bar b_j\}_{j = 1}^n$, write $b_{ij\bar k} := b_i \cdot b_j \cdot \bar b_k$, then $\Theta$ becomes
 \begin{equation*}
  \begin{split}
   \Theta = & \frac{1}{6} \sum_{i,j,k} \Lambda(b_i, b_j, \bar b_k)(b_{\bar i \bar j k} + b_{k \bar i \bar j} - b_{\bar i k \bar j}) + \frac{1}{6} \sum_{i,j,k} \Lambda(\bar b_i, \bar b_j, b_k)(b_{i j \bar k} + b_{\bar k i j}  - b_{i \bar k j})
  \end{split}
 \end{equation*}
 The terms in the first sum annilate $u_\mf l$. The terms in the second sum give
 \begin{equation*}
  \begin{split}
   & \left[b_{ij\bar k}, u_\mf l\right] = -(-1)^n u_\mf l \cdot b_i \cdot b_j\cdot \bar b_k  
   =  (-1)^n 2u_\mf l \cdot (\delta_{ik} b_j-\delta_{jk} b_i) \\
   & \left[b_{\bar k i j}, u_\mf l\right] = \bar b_k \cdot b_i \cdot b_j\cdot u_\mf l
   = 2(\delta_{ik} b_j - \delta_{jk} b_i) \cdot u_\mf l\\
   & \left[b_{i \bar k j}, u_\mf l\right] = b_i \cdot \bar b_k \cdot b_j \cdot u_\mf l - (-1)^n u_\mf l \cdot b_i \cdot \bar b_k \cdot b_j 
   =  2\delta_{jk} b_i \cdot u_\mf l - (-1)^n 2 u_\mf l \cdot \delta_{ik} b_j
  \end{split}
 \end{equation*}
 It follows that
 \begin{equation*}
  \begin{split}
   & d^\Cl u_\mf l = \frac{1}{4}\left[\Theta, u_\mf l\right] = \frac{1}{4} \sum_{i,j}\Lambda(\bar b_i, \bar b_j, b_i) ( b_j \cdot u_\mf l + (-1)^n u_\mf l \cdot b_j) \\
   = & \frac{1}{4} \sum_{i,j} \kappa([b_i, \bar b_i], \bar b_j)  (b_j \cdot u_\mf l + (-1)^n u_\mf l \cdot b_j)  
   = \frac{1}{2}\left(\rho \cdot u_\mf l + (-1)^n u_\mf l \cdot \rho\right)
  \end{split}
 \end{equation*}
 The last expression coincides with the Clifford action of $\displaystyle{\frac{1}{2}(-\rho, \rho) \in \Cl(\bar{\mf d})}$ on $u_\mf l$.
 \qed

\begin{corollary}\label{coro:ClifforddifferentialSamelsonandconjugate}
 Write $\mf l := \mf l_+$, then
 $\displaystyle{d^\Cl(u_\mf l \cdot u_\bmf l) = \frac{1}{2}(\rho, \rho) \circ (u_\mf l \cdot u_\bmf l)}$.
\end{corollary}
{\it Proof:}
 By Lemma \ref{lemma:ClifforddifferentialSamelson}, we have $\displaystyle{d^\Cl u_\bmf l = \frac{1}{2}(\rho, -\rho) \circ u_\bmf l}$. The identity follows.
\qed

We have the following more general statement.
\begin{prop}\label{prop:ClifforddifferentialtwoSamelson}
 Let $\mf l_\pm \subset \mf b_\pm$ be a pair of Samelson algebras and $\rho_\pm \in \mf b_\pm$ be the respective Weyl vectors. Let $u_+ \in \Cl(\mf g, \kappa)$ be a spinor defining $\bmf l_- \bdsum \bmf l_+$, then $\displaystyle{d^\Cl u_+ = \frac{1}{2}(-\rho_-, \rho_+) \circ u_+}$.
\end{prop}
{\it Proof:} First suppose that $\mf b_\pm = \mf b$, then we have $\bmf m = \bmf l_+ \dsum \bmf p$ with $\bmf p \subseteq \mf h \subseteq \mf b$.
Let $h \in \mf h$,
$$d^\Cl h = - \frac{1}{4}\sum_{\alpha \in R_+} \left( [h, a_\alpha] \cdot a_{-\alpha} + [h, a_{-\alpha}] \cdot a_\alpha\right) = -\frac{1}{4} \sum_{\alpha \in R_+} \alpha(h)(a_{\alpha} \cdot a_{-\alpha} - a_{-\alpha} \cdot a_{\alpha})$$
Suppose that $h \in \mf h$ and $h \not \in \mf l$. Write $\mf l := \mf l_+$, then
\begin{equation*}
 \begin{split}
  & d^\Cl(u_\mf l \cdot h) = d^\Cl u_\mf l \cdot h + (-1)^n u_\mf l \cdot d^\Cl h \\
  = & \frac{1}{2}\rho \cdot u_\mf l \cdot h  + (-1)^n \frac{1}{2}u_\mf l \cdot \rho \cdot h + (-1)^n u_\mf l \cdot \left(-\frac{1}{4}\right) \sum_{\alpha \in R_+} \alpha(h)(2 - 2a_{-\alpha} \cdot a_{\alpha}) \\
  = &  \frac{1}{2}\rho \cdot u_\mf l \cdot h  + (-1)^{n+1} u_\mf l \cdot \frac{1}{4}\sum_{\alpha \in R_+}(2\kappa(h_\alpha, h) - h_\alpha \cdot h) \\
  = & \frac{1}{2}\rho \cdot u_\mf l \cdot h + (-1)^{n+1} \frac{1}{2}u_\mf l \cdot h \cdot \rho =\frac{1}{2} (-\rho, \rho) \circ (u_\mf l \cdot h)
 \end{split}
\end{equation*}
The result then follows from mathematical induction on the dimension of $\mf p$.

In general, let $\hat g\in D$ such that $\Ad_{\hat g}(\mf b_- \bdsum \mf b_+) = \mf b \bdsum \mf b$, then $\hat g \circ u_+$ is a spinor defining $\Ad_{\hat g}(\bmf l_- \bdsum \bmf l_+)$. From previous computation, we have
$$\hat g \circ d^\Cl u_+ = d^\Cl(\hat g \circ u_+) = \frac{1}{2}(-\rho, \rho) \circ (\hat g \circ u_+)$$
Thus
$$d^\Cl u_+ = \frac{1}{2}\Ad^{-1}_{\hat g} (-\rho, \rho) \circ u_+$$
and the statement follows.
\qed

\begin{theorem}\label{thm:canonicalbundledegrees}
 Let $(K, \gamma; \JJ_\pm, \sigma)$ be a Lie algebraic generalized K\"ahler structure, and $U_+$ be the $\JJ_+$-canonical line bundle. Then 
$$\deg_\pm(U_+) = -2|\rho|^2$$
\end{theorem}
{\it Proof:} 
 The degrees are well defined according to Lemma \ref{lemma:alwaysGauduchon}.
 We prove it for induced generalized K\"ahler structures, and the general case is identical, with more cumbersome notations concerning the respective Weyl vectors.
 Furthermore, we consider only the $+$-degree, and the $-$-degree is similar. 
 Let $\chi_+ = \mc Q(u_+)$, by Lemata \ref{lemma:spinortrivialization} and \ref{lemma:differentialtrivialization}, we have
 $$d_\gamma \chi_+ = \frac{1}{2} \mf X_{(-\rho, \rho)} \cdot \chi_+$$
 Let $D_{U_+}$ denote the structure of $\JJ_+$-holomorphic line bundle on $U_+ = \<\chi_+\>$. For $s = \mf X_{(a, a')}$ with $a \in \bmf l_-$ and $a' \in \bmf l_+$, we have
 $$D_{U_+, s} \chi_+ = \mf X_{(a, a')} \cdot \mf X_{(-\rho, \rho)} \cdot \chi_+ = \kappa(a'+a, \rho) \chi_+$$
 Let $\rho^{1,0}_\pm \in \mf l_\pm$ denote the projections of $\rho$ to $\mf l_\pm$ respectively, and
 \begin{equation}\label{eq:canonicalbundleconnectionform}
  \varphi_+ = \kappa(\theta^r, \rho^{1,0}_+) \text{ and } \varphi_- = -\kappa(\theta^l, \rho^{1,0}_-)  
 \end{equation}
 Then, the induced structures of $I_\pm$-holomorphic line bundles are
 \begin{equation}\label{eq:canonicalbundlebiholomorphicstruct}
  D_+ \chi_+ = \varphi_+ \tensor \chi_+ \text{ and } D_- \chi_+ = \varphi_- \tensor \chi_+
 \end{equation}
 
 A natural Hermitian metric on $U_+$ is defined by $|\chi_+|^2 = 1$. The Chern connection for $D_+$ is
  $$\Cnabla{+} \chi_+ = (\varphi_+ - \bar\varphi_+) \tensor \chi_+$$
 whose curvature is given by
 $$F_+ = d(\varphi_+ - \bar \varphi_+)$$
 Use the basis $\{b_j\}_{j=1}^n$ of $\mf l := \mf l_+$ and the dual basis $\{\bar b_j\}$ of $\bmf l$. Then $\{b_j, \bar b_j\}$ gives a basis of $\mf g$. Let $\{b_j^*, \bar b_j^*\}$ be the dual basis of $\mf g^*$, e.g. $\<b_i^*, b_j\> = \delta_{ij}$. Then the value of $\omega_+$ at $e \in K$ given by
 $$\omega_+(e) = \frac{i}{2} \sum_{i} b_i^* \wedge \bar b_i^*$$
 To compute the degree $\deg_+(U_+)$, we are interested in $F_+ \wedge \omega_+^{n-1}$.
 Since all objects involved are right invariant, we only have to carry out the computation at $e \in K$, i.e. on the level of Lie algebra.
 The $b^*_i \wedge \bar b_i^*$ terms of $d\varphi_+(e)$ are
 \begin{equation*}
  \begin{split}
   & d\kappa(\theta^r, \rho^{1,0}_+) = -\kappa([\theta^r, \theta^r], \rho^{1,0}_+) \\
   \Longrightarrow & -\sum_{i, j}\kappa([b_i, \bar b_i], \rho^{1,0}_+) b^*_i \wedge \bar b^*_i = -\sum_{\alpha \in R_+} \kappa(h_\alpha, \rho^{1,0}_+) b_\alpha^* \wedge \bar b_\alpha^* \\
   = & -\sum_{\alpha \in R_+} \kappa(h_{\alpha,+}^{0,1}, \rho^{1,0}_+) b_\alpha^* \wedge \bar b_\alpha^* 
  \end{split}
 \end{equation*}
 It follows that the $b^*_i \wedge \bar b_i^*$ terms of $F_+(e)$ are
 \begin{equation*}
   -\sum_{\alpha \in R_+} \left(\kappa(h_{\alpha,+}^{0,1}, \rho^{1,0}_+)  + \bar{\kappa(h_{\alpha,+}^{0,1}, \rho^{1,0}_+)} \right) b_\alpha^* \wedge \bar b_\alpha^* = -2\sum_{\alpha \in R_+}\Re(\kappa(h_{\alpha,+}^{0,1}, \rho^{1,0}_+))b_\alpha^* \wedge \bar b_\alpha^*
 \end{equation*}
 and
 $$F_+(e) \wedge \omega_+^{n-1}(e) = -\frac{2 n! i^{n-1}}{2^{n-1}}\sum_{\alpha \in R_+} \Re(\kappa(h_{\alpha,+}^{0,1}, \rho^{1,0}_+)) \mathop\wedge_i (b_i^* \wedge \bar b_i^*) = -\frac{4n!i^{n-1}}{2^{n-1}}\Re(\kappa(\rho_+^{0,1}, \rho_+^{1,0})) \mathop\wedge_i (b_i^* \wedge \bar b_i^*)$$
 Since $\kappa(\rho, \rho) = 2 \Re\kappa(\rho_+^{0,1}, \rho_+^{1,0})$ we have
 $$\deg_+(U_+) = \frac{i}{2}\frac{\displaystyle{\int_K F_+ \wedge \omega_+^{n-1}}}{\displaystyle{\int_K \omega_+^n}} = \frac{i}{2}\frac{\displaystyle{-\frac{2n!i^{n-1}}{2^{n-1}}\kappa(\rho, \rho) \mathop\wedge_i (b_i^* \wedge \bar b_i^*)}}{\displaystyle{\frac{n!i^{n}}{2^{n}} \mathop\wedge_i (b_i^* \wedge \bar b_i^*)}} = -2\kappa(\rho, \rho)$$
\qed

\begin{remark}\label{remark:positiveproperty}
 \rm{
  Although Proposition \ref{prop:ClifforddifferentialtwoSamelson} and Theorem \ref{thm:canonicalbundledegrees} concern only the structure $\JJ_+$, it's easy to see that similar results hold for $\JJ_-$. Namely, we have
  $$\deg_\pm(U_-) = - 2|\rho|^2$$
  where the $\deg_\pm$ in this case are respect to $\pm I_\pm$.
  These results may be interpreted as a certain \emph{positivity} of the Lie algebraic generalized K\"ahler structures.
 }
\end{remark}

We recover the following result from \cite{Cavalcanti12}.
\begin{corollary}\label{coro:canonicalgenKahlerisCalabiYau}
 Let $(K, \gamma; \JJ_\pm, \sigma)$ be a Lie algebraic generalized K\"ahler structure. Then $\JJ_+$ is generalized Calabi-Yau iff $K$ is abelian, i.e. a torus. 
\end{corollary}
{\it Proof:}
 Since $\rho = 0 \iff K$ is a torus, we see that the $\JJ_+$-canonical bundle is holomorphically trivial iff $K$ is a torus.
\qed

\begin{corollary}\label{coro:torusgenCalabiYau}
 Let $T$ be an even dimensional torus, then any invariant generalized K\"ahler structure on $T$ is a generalized Calabi-Yau metric structure.
\end{corollary}
{\it Proof:}
 For a torus, $\rho = 0$ and $\chi_+$ is a non-vanishing holomorphic section of $U_+$, which implies that $\JJ_+$ is generalized Calabi-Yau. That $\JJ_-$ is also generalized Calabi-Yau in this case follows from Remark \ref{remark:positiveproperty}. Proposition $4.2$ $d)$ in \cite{AlekseevBursztynMeinrenken} gives the comparison of lengths of the spinors.
\qed

\subsection{Hodge decomposition}\label{subsect:CliffordHodgedecomposition}
The Hodge decomposition of the twisted de Rham complex $(\Omega^*(K), d_\gamma)$ %
induces the \emph{Hodge decomposition} of $\Cl(\mf g, \kappa)$. 
Consider the spinor actions of $\hat J_+ := J_-\dsum J_+$ and $\hat J_- := -J_-\dsum J_+$ on $\Cl(\mf g, \kappa)$.
\begin{lemma}\label{lemma:spinorJaction}
 Let $\{b_{j,\pm}\}$ be a complex basis of $\mf l_\pm$ and $\{\bar b_{j,\pm}\}$ the dual basis of $\bmf l_\pm$ with respect to $\kappa$, i.e. $\kappa(b_{i,\pm}, \bar b_{j,\pm}) = \delta_{ij}$. Define $\tau_{J_\pm} \in \Cl(\mf g, \kappa)$ by
 $$\tau_{J_+} = -\frac{ni}{2} + \frac{i}{2}\sum_{j= 1}^n b_{j,+}\cdot \bar b_{j,+} \text{ and } \tau_{J_-} = \frac{ni}{2} - \frac{i}{2}\sum_{j = 1}^n \bar b_{j,-} \cdot b_{j,-}$$
 Let $\hat\tau_\pm = (\pm\tau_{J_-}, \tau_{J_+})$ then the action of $\hat J_\pm$ on $\Cl(\mf g, \kappa)$ is given by
 $$u \mapsto \hat \tau_\pm \circ u = \tau_{J_+}\cdot u \mp u \cdot \tau_{J_-}$$
\end{lemma}
{\it Proof:}
 Since $\kappa$ is Hermitian with respect to both $J_\pm$, $\tau_{J_\pm}$ are the corresponding elements in $\Cl(\mf g, \kappa)$ representing $J_\pm$:
 $$\tau_{J_\pm} = \frac{1}{4}\sum_{j = 1}^n \left(J_\pm\left(b_{j,\pm}\right)\cdot \bar b_{j,\pm} + J_\pm\left(\bar b_{j,\pm}\right)\cdot b_{j,\pm}\right)$$
 The statements follow from straightforward computations and the definition of spinor actions.
\qed

\begin{corollary}\label{coro:canonicalspinoraction}
 Let $u_+$ be the pure spinor defining $\bmf l_- \bdsum \bmf l_+$, then $\hat \tau_+ \circ u_+ = -in u_+$ and $\hat \tau_- \circ u_+ = 0$.
\end{corollary}
{\it Proof:} 
 It follows from \eqref{eq:pluspinordefn} that $\displaystyle{\tau_{J_+} \cdot u_+ = -\frac{ni}{2} u_+}$ and $\displaystyle{u_+ \cdot \tau_{J_-} = \frac{ni}{2} u_+}$.
\qed

Since $\mf l_\pm$ are isotropic with respect to $\kappa$, we have $\wedge^*\mf l_\pm \subset \Cl(\mf g, \kappa)$.
\begin{prop}\label{prop:CliffordHodgedecomposition}
 $\displaystyle{\Cl(\mf g, \kappa) = \dsum_{r,s} \mf U_{r,s}}$, where $\mf U_{r,s} := \wedge^p \mf l_+ \cdot u_+ \cdot \wedge^q\mf l_-$ with $r = p+q - n$ and $s = p-q$. Moreover, $\mf U_{r,s}$ is the $i(r,s)$-eigensubspace of the spinor action by $(\hat\tau_+, \hat\tau_-)$.
\end{prop}
{\it Proof:}
 First we note that $\displaystyle{\Cl(\mf g, \kappa) = \Span\{\mf U_{r,s}\}}$. Then the second statement implies the first one.
 For any $p,q = 0, 1, \ldots, n$, we compute the action of $\hat\tau_+$ on 
  $$u := v_+ \cdot u_+ \cdot v_- := b_{1,+}\cdot\ldots\cdot b_{p,+} \cdot u_+ \cdot b_{1,-} \cdot\ldots\cdot b_{q,-}$$
 and show that it is an $ir$-eigenvector, with $r = p+q - n$. The cases where $v_\pm$ are general basis elements of $\wedge^* \mf l_\pm$ are similar.
 \begin{equation*}
  \begin{split}
   & \hat\tau_+\circ u = \tau_{J_+}\cdot u - u \cdot \tau_{J_-} = -in u + \frac{i}{2}\sum_{j=1}^n \left(b_{j,+}\cdot \bar b_{j,+} \cdot u + u \cdot \bar b_{j,-}\cdot b_{j,-}\right)\\
   = & -in u + \frac{i}{2}\sum_{j=1}^n \left(b_{j,+}\cdot \bar b_{j,+} \cdot v_+\cdot u_+\cdot v_- + v_+\cdot u_+\cdot v_- \cdot \bar b_{j,-}\cdot b_{j,-}\right) \\
   = & -in u + \frac{i}{2}\sum_{j=1}^p(2v_+\cdot u_+ \cdot v_-) + \frac{i}{2}\sum_{j=1}^q(2v_+\cdot u_+ \cdot v_-)\\
   = & i(p+q - n) u = iru
  \end{split}
 \end{equation*}
 The computation for the action of $\hat\tau_-$ is similar.
\qed

\subsection{Cohomology of $\bar L$} \label{subsubsect:cohomologyofLonK}
From now on, 
we write $\JJ := \JJ_+$, $L: = L_+$.

Let $\{\xi_i^\pm : i = 1, \ldots, n\}$ be a basis of the dual $\bar{\mf l}_\pm^*$ of the complex Lie subalgebra $\bar{\mf l}_\pm \subset \mf g$, respecting the decompositions $\bar{\mf l}_\pm = \mf t_{\pm, (0,1)} \dsum \mf r_\pm$. For $i = 1, \ldots, n$, let $\alpha_i$ (respectively, $\beta_i$) be the right invariant (respectively, left invariant) $1$-form on $K$ such that $\alpha_i(e) = \xi^+_i$ (respectively, $\beta_i(e) = \xi^-_i$). Then $\{\alpha_i\}_{i=1}^n$ is the basis of right invariant $I_+$-$(0,1)$-forms and $\{\beta_i\}_{i = 1}^n$ is the basis of left invariant $I_-$-$(0,1)$-forms. 
We use the same notations, $\alpha_i$ and $\beta_i$, to denote the corresponding elements in $\Omega^1(\bar\ell_\pm)$. For $P \subset \{1, 2, \ldots, n\}$, we use $\xi^+_P$ and $\beta_P$ to denote $\displaystyle{\mathop\wedge_{i \in P} \xi^+_i}$ and $\displaystyle{\mathop\wedge_{i \in P} \beta_i}$, and etc.
\begin{prop}\label{prop:Liealgebroidastrivialholobundle}
 $\mc A_\pm$ are trivial as $I_\pm$-holomorphic vector bundles.
\end{prop}
{\it Proof:}
 By definition, $\{\alpha_i\}$ and $\{\beta_i\}$ are global frames of $\bar \ell_+^*$ and $\bar \ell_-^*$ respectively. By Lemma \ref{lemma:Bismutflat}, $\bar\delta_- \alpha_i = 0$ and $\bar\delta_+ \beta_i = 0$ for all $i$. Lemma \ref{lemma:explicitreducedholomorphicstructure} implies that they are global holomorphic frames of $\mc A_\pm^*$ respectively. Thus $\mc A_\pm$ are holomorphically trivial.
\qed

\begin{corollary}\label{coro:firstpageforLonK}
 The first page of the first spectral sequence associated to $\Omega^{p,q}(\bar L)$ is
 $${}_IE^{p,q}_1 = \left(\wedge^p \bar{\mf l}_+^* \tensor \wedge^q \mf t_-^{0,1}, d_{\bar{\mf l}_+} \tensor \id \right)$$
 where $\mf t_-^{0,1} = \left(\mf t_{-,(0,1)}\right)^*$.
\end{corollary}
{\it Proof:} The description of the first page of the spectral sequence in \S \ref{subsect:holomorphicreduction} gives
$${}_IE^{p,q}_1 = H^{0,q}_{\bar\partial_-}(\wedge^q \mc A_-^*)$$
By Proposition \ref{prop:Liealgebroidastrivialholobundle}, $\wedge^q \mc A_-$ is trivial as well. The explicit trivialization gives
$${}_IE^{p,q}_1 = H^{0,q}_{\bar\partial_-}(\wedge^q \mc A_-^*) = \wedge^p \bar{\mf l}_+^* \tensor H^{0,q}_{\bar\partial_-}(K)$$
The expression of the $E_1$-terms then follows from the general fact (cf. \cite{Pittie88}) that
$$H^{0,q}_{\bar\partial_-}(K) \cong \wedge^q \mf t_-^{0,1}$$

To compute the induced differential $d_1$, we note that $\bar\partial_+$ coincides with the differential $d_{\bar{\mf l}_+}$ on $\wedge^p \bar{\mf l}_+^*$. By Alexandrov-Ivanov \cite{AlexandrovIvanov}, $H^{0,q}_{\bar\partial_-}(K)$ can be represented by $\bar\partial_-$-Harmonic forms, which are left invariant. From Corollary \ref{coro:formBismutflat}, $\bar\partial^- = 0$ on $\wedge^q \mf t_-^{0,1}$. Thus, $\bar\delta_+$ induces $d_1 = d_{\bar{\mf l}_+} \tensor \id$.
\qed 

\begin{prop}\label{prop:secondpageforLonK}
 The first spectral sequence associated to $\Omega^{p,q}(\bar L)$ collapses at $E_2$, where
 $${}_IE_2^{p.q} = \wedge^p \mf t_+^{0,1} \tensor \wedge^q \mf t_-^{0,1}$$
\end{prop}
{\it Proof:}
 Either using the Hochschild-Serre spectral sequence (Hochschild-Serre \cite{HochschildSerre53}) or quoting \cite{Pittie88} again, we have $H^*(\bar{\mf l}_+) \cong H^*(\mf t_+^{0,1}) = \wedge^*\mf t_+^{0,1}$. To obtain the differential $d_2$, we note that a basis of ${}_IE_2^{p.q}$ is given by the classes represented by the forms $\alpha_Q \wedge \beta_{Q'}$, where $|Q| = p$, $|Q'| = q$, such that $\xi^+_Q \in \wedge^p \mf t_+^{0,1}$ and $\xi^-_{Q'} \in \wedge^q \mf t_-^{0,1}$. Since
 $$\bar\delta_+(\alpha_Q \wedge \beta_{Q'}) = \bar\partial_+ \alpha_Q \wedge \beta_{Q'} + (-1)^p \alpha_Q \wedge \bar\partial^+ \beta_{Q'} = 0$$
 we see that $d_2 = 0$ on ${}_IE_2^{p.q}$.
\qed

From Proposition \ref{prop:isomorphiccohomology}, we obtain:
\begin{corollary}\label{coro:cohomologyofLonK}
 $\HH^*(\mc A_+) \cong \HH^*(\mc A_-) \cong H^*(\bar L) \cong \wedge^*(\mf t_-^{0,1} \bdsum \mf t_+^{0,1})$.
 \qed
\end{corollary}

\subsection{$\JJ$-Picard group}\label{subsubsect:JplusPicard}
Let $V = K \times \CC$ be the topologically trivial complex line bundle on $K$. The $I_+$-holomorphic structures on $V$ is then classified by the classical Picard group with respect to $I_+$
$$\Pic^+_0(K) = \frac{H^1(K; \mc O_+)}{H^1(K; \ZZ)}$$
By \cite{Pittie88}, $H^1(K; \mc O_+) \cong \CC^r$ where $2r$ is the rank of $K$. In general, $H^1(K; \ZZ)$ may not be a full lattice in $H^1(K; \mc O_+)$. For example, when $K$ is semi-simple and simply connected, $H^1(K; \ZZ) = \{0\}$ and $\Pic^+_0(K) \cong \CC^r$. In general, when $H^1(K;\ZZ)$ has no torsion, the Picard group is of the form
$$\Pic^+_0(K) \cong \CC^{r_1} \times \left(\CC^\times\right)^{r_2} \times (S^1)^{2r_3}, \text{ with } r_1 + r_2 + r_3 = r$$
For example, let $T = (S^1)^2$ be the $2$-torus, $H = SU(2) \times S^1$ the Hopf surface, then we have
$$\Pic_0(T) \cong (S^1)^2, \Pic_0(H) \cong \CC^\times \text{ and } \Pic_0(SU(3)) \cong \CC$$

Let $\PPic_0(K)$ denote the space classifying $\JJ$-holomorphic structures on $V$. Then $\PPic_0(K)$ is a group under the tensor product, with identity represented by the trivial connection 
$$D_0: C^\infty(V) \to C^\infty(V \tensor \bar L^*) : D_0(u) = 0$$
where $u$ is the constant section with value $1$. We note that in this case, $C^\infty(V \tensor \bar L^*) \cong \Omega^1(\bar L)$ and $D_0$ coincides with $d_{\bar L}$. Any $\bar L$-connection $D$ on $V$ is of the form
$$D = d_{\bar L} + \omega, \text{ with } \omega \in \Omega^1(\End(V); \bar L) \cong \Omega^1(\bar L)$$
The flatness of $D$ is equivalent to the Maurier-Cartan equation
$$d_{\bar L}\omega + \omega \wedge \omega = d_{\bar L} \omega = 0$$
It then follows that
$$\PPic_0(K) \cong \frac{\ker d_{\bar L}}{\mc G(V)}$$
where $\mc G(V) = C^\infty(\Aut(V))$ is the gauge group acting on the space of $\bar L$-connections:
$$g\circ D := D + g^{-1} d_{\bar L} g$$
Infinitesimally, the tangent space of $\PPic_0(K)$ at $D_0$ is given by $H^1(\bar L) \cong \mf t_+^{0,1} \dsum \mf t_-^{0,1} \cong \CC^{2r}$.

On any generalized K\"ahler manifold, Lemma \ref{lemma:Picardhomomorphism}, together with the isomorphism of categories in Proposition \ref{prop:isomorphiccategories}, gives the exact sequence of complex analytic group homomorphisms:
$$0 \to \Pic^{\mc A_+}_0 \to \PPic_0 \to \Pic^+_0$$
On $K$, the sequence above can be completed to a split short exact sequence.
\begin{prop}\label{prop:JPicardonK}
 $\PPic_0(K) \cong \Pic^+_0(K) \times \Pic^{\mc A_+}_0(K) \cong \Pic^+_0(K) \times \CC^r$.
\end{prop}
{\it Proof:}
 We only have to show that $\Pic^{\mc A_+}_0(K) \cong \CC^r$, then dimension counting implies that the forgetful map $\PPic_0(K) \to \Pic^+_0(K)$ is a submersion. The result then follows.
 
 On the trivial $I_+$-holomorphic bundle $\mc O_+$, let $\partial_{\mc A_+, 0}$ be the trivial $\mc A_+$-connection. Then any flat $\mc A_+$-connection on $\mc O_+$ is of the form
 $$\partial_{\mc A_+} = \partial_{\mc A_+, 0} + \alpha, \text{ with } \alpha \in \mc A_+^* \text{ and } \partial_{\mc A_+,0} \alpha = 0$$
 We note that the only holomorphic bundle automorphisms on $\mc O_+$ are multiplications by elements of $\CC^\times$. Since $\mc A_+$ is a trivial $I_+$-holomorphic bundle by Proposition \ref{prop:Liealgebroidastrivialholobundle}, similar arguments as in the proof of Proposition \ref{prop:secondpageforLonK} gives $\Pic^{\mc A_+}_0(K) \cong \CC^r$. 
\qed

\begin{remark}\label{remark:JPicardfromBH}
 \rm{
  The $\JJ$-Picard group on a generalized K\"ahler manifold can also be seen from the bi-Hermitian interpretation as following. By Lemma \ref{lemma:commutationrelationdefinesgenholo}, $V$ is a $\JJ$-holomorphic line bundle iff it is $I_\pm$-holomorphic and the induced differentials $D_\pm$ satisfy the commutation relation. For $V = M \times \CC$, by an abuse of notations, let $D_0$ denote both of the trivial $I_\pm$-holomorphic structures on $V$, then an $I_\pm$-holomorphic structure is given by a pair
  $$(\alpha_+, \alpha_-) \in \Omega^1_+(M) \times \Omega^1_-(M) \text{ such that } \bar\partial_+\alpha_+ = \bar\partial_-\alpha_- = 0$$
  where the complex gauge group $\mc G$ acts diagonally.
  The commutation relation in this case becomes
  $$\bar\delta_-\alpha_+ + \bar\delta_+\alpha_- = 0$$
  Thus, the $\JJ$-Picard group can also be seen as the quotient
  $$\PPic_0(M) = \frac{\{(\alpha_+, \alpha_-) : \bar\delta_-\alpha_+ + \bar\delta_+\alpha_- = 0\}}{\mc G}$$
  This is in general different from $\Pic^+_0(M) \times \Pic^-_0(M)$.
 }
\end{remark}


\begin{thebibliography}{10}
\bibitem{AlekseevBursztynMeinrenken} A.~Alekxeev, H.~Bursztyn and E.~Meinrenken, Pure spinors on Lie groups, {\it Ast\'erisque} No. 327 (2009), pp. 131 -- 199. 
\bibitem{AlekseevskyDavid10} D.~V.~Alekseevsky and L.~David, Invariant generalized complex structures on Lie groups, {\it Proc. London Math. Soc.} (4) 105 (2012), pp. 703--729 
\bibitem{AlexandrovIvanov} B.~Alexandrov and S.~Ivanov, Vanishing theorems on Hermitian manifolds, {\it Differential Geom. Appl.} (3) 14 (2001), pp. 251--265.
\bibitem{Atiyah57} M.~F.~Atiyah, Complex analytic connections in fibre bundles, {\it Trans. Am. Math. Soc.}, 85 (1957), 181 -- 207.
\bibitem{Baraglia12} D.~Baraglia, Variation of Hodge structure for generalized complex manifolds, {\it Differential Geom. Appl.} 36 (2014), pp. 98--133 
\bibitem{Cavalcanti11} G.~Cavalcanti, The decomposition of forms and cohomology of complex manifolds of generalized complex manifolds, {\it J. Geom. and Phys.}
(1) 57 (2006), pp. 121--132.
\bibitem{Cavalcanti12} G.~Cavalcanti, Hodge theory and deformations of SKT manifolds, {\it arXiv:1203.0493v5}
\bibitem{EvensLu06} S.~Evens, J.-H.~Lu, On the variety of Lagrangian subalgebras, II, {\it Ann. Sci. cole Norm. Sup.} (4) 39 (2006), no. 2, 347 -- 379.
\bibitem{Ferreira13} A.~C.~Ferreira, A vanishing theorem in twisted De Rham cohomology, {\it Proceedings of the Edinburgh Mathematical Society (Series 2)}, (2) 56 (2013), pp. 501 -- 508.
\bibitem{GatesHullRocek} S.~J.~Gates, Jr., C.~M.~Hull, and M.~Ro\v cek, Twisted multiplets and new supersymmetric nonlinear $\sigma$-models, {\it Nuclear Phys. B} (1) 248 (1984), pp. 157 -- 186.
\bibitem{Gualtieri04} M.~Gualtieri, Generalized complex geometry, Ph.D. thesis, Oxford, November 2003 {\it arXiv:math/0401221v1}
\bibitem{Gualtieri0409} M.~Gualtieri, Generalized geometry and the Hodge decomposition, {\it arXiv:math/0409093v1}
\bibitem{Gualtieri0703} M.~Gualtieri, Generalized complex geometry, {\it Ann. Math.} (1) 174 (2011) pp. 75 -- 123. 
\bibitem{Gualtieri0710} M.~Gualtieri, Branes on Poisson varieties, {\it The Many Facets of Geometry} 1 (9), pp. 368--395. 
\bibitem{Gualtieri10} M.~Gualtieri, Generalized K\"ahler geometry, {\it Comm. Math. Phys.} (1) 331 (2014), pp. 297--331. 
\bibitem{Hitchin02} N.~Hitchin, Generalized Calabi-Yau manifolds,  {\it Quart. J. Math. Oxford} 54 (2003), pp. 281--308.
\bibitem{Hitchin06} Nigel~Hitchin, Brackets, forms and invariant functionals, {\it Asian J. Math.}, 10 (2006), no. 3, 541 -- 560. 
\bibitem{HochschildSerre53} G.~Hochschild and J.~P.~Serre, Cohomology of Lie algebras, {\it Ann. Math.}, 57 (1953), 591 -- 603.
\bibitem{HuMoraruSeyyedali} S.~Hu, R.~Moraru and R.~Seyyedali, A Kobayashi-Hitchin correspondence for $I_\pm$-holomorphic bundles, {\it  arXiv:1411.3416}.
\bibitem{HuUribe06} S.~Hu and B.~Uribe, Extended manifolds and extended equivariant cohomology, {\it J. Geom. Phys.} 59 (2009), pp. 104 -- 131. 
\bibitem{Humphreys} J.~E.~Humphreys, Introduction to Lie algebras and representation theory, {\it Graduate texts in mathematics} 9, Springer-Verlag, 1972.
\bibitem{Karolinsky00} E.~Karolinsky, A classification of Poisson homogeneous spaces of complex reductive Poisson-Lie groups, {\it Banach Center Publ.} 51. Polish Acad. Sci., Warsaw, 2000.
\bibitem{LaurentGengouxStienonXu08} C.~Laurent-Gengoux, M.~Sti\'enon, and P.~Xu. Holomorphic Poisson manifolds and holomorphic Lie algebroids. {\it Int. Math. Res. Not. IMRN}, Art. ID rnn 088, 46, 2008.
\bibitem{Lu97} J.-H.~Lu, Poisson homogeneous spaces and Lie algebroids associated to Poisson actions, {\it Duke Math. J.}, 86(1997), no. 2, 261--304.
\bibitem{Mackenzie07} K.~C.~H.~Mackenzie, Ehresmann doubles and Drinfel’d doubles for Lie algebroids and Lie bialgebroids, {\it Journal f\"ur die reine und angewandte Mathematik (Crelles Journal)}, 658 (2011), pp. 193--245. 
\bibitem{Meinrenken13} E.~Meinrenken, Clifford algebras and Lie theory, {\it Ergebnisse der Mathematik und ihrer Grezgebiet. 3. Folge, Vol. 58}, Springer-Verlag, 2013.
\bibitem{Mokri97} T.~Mokri, Matched pairs of Lie algebroids, {\it Glasgow Math. J.}, (2) 39 (1997), pp. 167 -- 181.
\bibitem{Pittie88} H.~R.~Pittie, The Dolbeault-cohomology ring of a compact, even-dimensional lie group, {\it Proc. Indian Acad. Sci. (Math. Sci.)}, 98 (1988), pp. 117--152.
\bibitem{Samelson53} H.~Samelson, A class of complex analytic manifolds, {\it Portugaliae Math.}, 12 (1953), pp. 129 -- 132.
\bibitem{Tortella11} P.~Tortella, $\Lambda$-modules and holomorphic Lie algebroid connections, {\it Centr. Eur. J. Math.} 10 (4), pp. 1422--1441 
\bibitem{Wang54} H.~C.~Wang, Closed manifolds with homogeneous complex structure, {\it Amer. J. Math.}, 76 (1954), pp. 1--32.
\bibitem{Wells79} R.~O.~Wells,~Jr., Differential analysis on complex manifolds, {\it Graduate texts in mathematics} 65, Springer-Verlag, 1979.
\end{thebibliography}
\end{document}